\documentclass[final,leqno]{siamltex}

\usepackage{color}

\newcommand{\DC}[1]{{\textcolor{black}{#1}}}

\usepackage{amssymb,amsmath,amsfonts}
\usepackage{epsfig}
\usepackage{psfrag}
\usepackage{enumerate}

\usepackage[pdftex,hidelinks]{hyperref}
\hypersetup{a4paper = true,
  pdfpagemode = None,
  pdfstartview  = FitH}

\title{A direction preserving discretization for computing phase-space densities}

\author{David J. Chappell\footnotemark[2], Jonathan J. Crofts\footnotemark[2], Martin Richter\footnotemark[2],\footnotemark[3] \and Gregor Tanner\footnotemark[3]\:\:\thanks{Support from the EPSRC (grant no.~EP/R012008/1) is
gratefully acknowledged. The authors also gratefully acknowledge the shock tower mesh data provided by Stephen Fisher, Jaguar Land Rover Ltd.}}

\begin{document}

\maketitle
\renewcommand{\thefootnote}{\fnsymbol{footnote}}

\footnotetext[2]{School of Science and Technology, Nottingham Trent
University, Clifton Campus, Clifton Lane,
        Nottingham, UK NG11 8NS ({\tt david.chappell@ntu.ac.uk}).}
\footnotetext[3]{School of Mathematical Sciences, University of
Nottingham, University Park, Nottingham NG7 2RD, UK.}
\begin{abstract}
Ray flow methods are an efficient tool to estimate vibro-acoustic or electromagnetic energy transport in complex domains at high-frequencies. Here, a Petrov-Galerkin discretization of a phase-space boundary integral equation for transporting wave energy densities on two-dimensional surfaces is proposed. The directional dependence of the energy density is approximated at each point on the boundary in terms of a finite local set of directions propagating into the domain. The direction of propagation can be preserved for transport across multi-component domains when the directions within the local set are inherited from a global direction set. The range of applicability \DC{and computational cost} of the method will be explored through a series of numerical experiments, including wave problems from both acoustics and elasticity in both single and multi-component domains. The domain geometries considered range from both regular and irregular polygons to curved surfaces, including a cast aluminium shock tower from a Range Rover car. 
\end{abstract}



\pagestyle{myheadings} \thispagestyle{plain} \markboth{CHAPPELL, CROFTS, RICHTER AND TANNER}{DIRECTION PRESERVING DISCRETIZATION IN PHASE-SPACE}

\section{Introduction}
\subsection{Motivation and contribution}
Dynamical energy analysis (DEA) is a framework for modelling the flow of high-frequency wave energy densities that has been developed over the last ten years \cite{CTLS13, TH19, GT09}. At the heart of the method is a linear integral operator based model for transporting phase-space densities along ray trajectories between intersections with the boundary of a domain or sub-domain. In recent years, the capability of the method has extended to large-scale problems from industry through an efficient implementation on finite element type meshes (in the position variable) \cite{CTLS13, TH19}.  However, the benchmarking of the method against simple toy examples has often proved difficult \cite{JB17, CTLS13}. The problem stems from the fact that the lowest order DEA approximations represent the momentum (or direction) dependence of the energy density as a constant function and stronger directivity is incorporated by increasing the order of this approximation in some sense, usually by expanding in terms of Legendre polynomials \cite{JB17}. However, for simple examples with regular geometries, the dependence on direction is typically described in terms of a finite number of Dirac delta distributions and hence the direction basis approximation in DEA proves rather inefficient. 

The main contribution of this study is to propose an alternative Petrov-Galerkin discretization that will be efficient for approximating densities with the type of strong directional dependence that is usually problematic for DEA. In particular, the directional dependence will be approximated using a basis of Dirac delta distributions to specify a discrete finite set of directions oriented into the domain from the boundary. The propagation direction will be preserved when transporting densities through multi-domains, such as through a sequence of mesh cells, provided that the finite direction set at each boundary position is taken as a subset of a global set; a density can then continue along a ray with the same global direction in different sub-domains. As a consequence, the proposed methodology naturally extends to transport densities through complex geometries approximated by finite element type meshes, such as those arising in industry.  

\subsection{Wider context and literature survey}
The linear operator behind the DEA method is a modified form of the Frobenius--Perron (FP) operator $\mathcal{L}^{\tau}$. The FP operator describes the evolution through time $\tau$ of a density $f$ along the solution trajectories of a dynamical system defined by a vector field $\mathbf{V}$ as follows:
\begin{equation}\label{eq:DS}
\dot{X}=\mathbf{V}(X).
\end{equation}
The solutions of (\ref{eq:DS}) define trajectories in phase-space of the form $X(\tau)=\varphi^{\tau}(X(0))$, where $\varphi^{\tau}$ is the associated flow map. The FP operator may then be written as a linear integral operator
\begin{equation}\label{eq:FPO}
\mathcal{L}^{\tau}[f](X)=\int\delta(X-\varphi^{\tau}(Y))f(Y,0)\,\mathrm{d}Y.
\end{equation}
Depending on the underlying dynamical system, the FP operator (\ref{eq:FPO}) can be used to model a variety of physical phenomena including tracking uncertainties in the dynamics of hypersonic flight \cite{PD11}, analyzing the conformation dynamics of molecules \cite{GF09} and modes of optical cavities \cite{JK16} as well as finding regions of minimal external transport in atmospheric science \cite{Froyland10}. The underlying dynamical system for the DEA approach is a classical Hamiltonian system describing the propagation of ray trajectories in the phase-space $X = ({\bf r}, {\bf p})$ with position and momentum variables $\bf r$, $\bf p$, respectively \cite{GT09}.

Discretization methods for FP operators have traditionally focused on 1D dynamical systems, for example the classical Ulam method \cite{Ulam1964} whereby phase-space is divided into cells and the cell-to-cell transition rates are estimated. More recently, Junge \emph{et al.} proposed wavelet and spectral collocation methods for the infinitesimal FP operator \cite{FJK11, JK09} opening up the possibility to treat higher dimensional systems, but typically in simple geometric settings such as unit cubes or tori. Dynamic mode decomposition (DMD) based approaches have been put forward for a related transfer operator (Koopman) in \cite{BMM12} and subsequently extended to the FP operator and compared to Ulam's method, as well as generalised Galerkin methods, in a recent review article \cite{KKS16}. Here, extended DMD is shown to be efficient when the eigenfunctions are smooth and the authors suggest that extended DMD would then be suitable for high-dimensional systems, such as those arising in molecular dynamics. A further survey of data-driven model reduction methods for the FP (and Koopman) operator is given in \cite{KKS18}, including extended DMD along with methods developed by the dynamical systems and molecular dynamics communities such as time-lagged independent component analysis and generalised Markov state models. The analysis of high-dimensional systems is proposed in a problem specific manner via an informed choice of dimensionality reduction, tensor decomposition and sparsification (including deep learning) methods. Several other schemes have recently been proposed and rigorously analysed including the finite volume method \cite{RN18} and the Galerkin method \cite{QH19,CW19}. The discretization methods used in DEA previously have been based on Galerkin projections onto Legendre polynomial or Fourier basis expansions \cite{JB17, GT09}. The first rigorous results concerning the convergence of DEA appeared recently \cite{JS2020} for the case of Fourier basis expansions in circular domains. 

Instead of transporting densities along trajectories in phase-space, there are many alternative high-frequency wave models based on directly tracking rays, beams or wave-fronts. The focus of these direct approaches tends to be on scattering, refractive lens or aperture problems in free-space, or absorbent cavities, where only a small number of reflections must be modelled - a review of this family of methods can be found in Ref.~\cite{BE03}. Example applications can be found in a wide range of high-frequency wave problems including in room acoustics \cite{HK09, LS15}, electromagnetics \cite{HL89}, as well as in seismology \cite{Cer01} and underwater acoustics \cite{IS01}. Recent developments in terms of ray-based methods for heterogeneous problems include ray approximations of the heterogeneous Green's function for medical imaging applications \cite{FR19}, and the fast Huygens sweeping method for high-frequency Helmholtz \cite{WL16} and Maxwell \cite{JQ16} wave problems.

The DEA method computes the stationary density $\rho$ accumulated in the long-time limit of the dissipative FP operator:
\begin{equation}\label{eq:FPOstat}
\rho(X)=\lim_{T\rightarrow\infty}\int_0^T \int w(Y,\tau)\delta(X-\varphi^{\tau}(Y))f(Y,0)\,\mathrm{d}Y\mathrm{d}\tau,
\end{equation}
given an initial density distribution $f(X,0)$. %
The dissipative factor $w$ has been included to facilitate convergence in the long-time limit. The stationary density $\rho$ can then be used to describe the wave energy distribution corresponding to the geometrical optics (GO) limit of frequency-domain wave problems in finite domains \cite{GT09}. In this case, directly tracking rays or wave-fronts can quickly become intractable since multiple reflections at boundaries often lead to complicated folding patterns of the associated level-surfaces and an exponential increase in the number of branches that need to be considered. Instead, simplified statistical models are typically set up under additional ergodicity and mixing assumptions on the ray dynamics \cite{Gradoni2014, RL95}. One of the simplest and most well-known approaches of this type is the so-called Statistical Energy Analysis (SEA), which has proved popular for applications in vibration and acoustics \cite{Lyon1969, RL95}. The principle advantage of SEA models is that they describe highly complex systems using relatively small systems of equations by subdividing a built-up structure into a set of subsystems and considering the (thermodynamic) energy flow rates, or coupling loss factors, between subsystems. However, setting up an SEA model requires considerable expertise since there are a number of limiting assumptions on its validity that must be considered with care - these limitations are widely reported in the SEA literature, see for example Refs.~\cite{Lebot2010, Lafont2013}. DEA relaxes some of the limiting assumptions of SEA by including directional dependence as well as allowing for local variability in the positional dependence \cite{GT09}. As such, DEA provides a higher resolution result than SEA whilst at the same time removing the remodelling effort.

\subsection{Structure of the paper}

The paper is structured as follows:  in Section \ref{sec:PropDenOp}  we will describe a boundary integral operator model for transporting wave energy densities through phase-space, including the projection of the resulting boundary density into the domain.  Section \ref{sec:FiniteDimAppr} will then introduce the direction preserving Petrov-Galerkin discretization of both the boundary integral operator and the domain projection formula, whereby the latter can be reduced to a finite sum over the global direction set. The effectiveness of the proposed discretization scheme will then be demonstrated through the numerical experiments reported in Section \ref{sec:NumResults}. The numerical examples considered range from a simple problem in a unit square domain, through irregular coupled acoustic cavities and culminate in a simulation of the vibrational energy distribution within a triangle mesh of a thin aluminium shell taken from a vehicle shock tower.


\section{Propagating phase-space densities via integral operators}\label{sec:PropDenOp}

In this section we outline the underlying model for propagating phase-space densities
along ray trajectories using the Perron-Frobenius operator. We are concerned with transporting densities through multi-domains $\Omega=\cup_{j=1}^{K}\Omega_j$, via Hamiltonian ray dynamics $H_j(\mathbf{r},\mathbf{p})=|\mathbf{p}|/\eta(\mathbf{r})\equiv1$, where $\eta(\mathbf{r})>0$ is related to the inverse of the phase velocity (or slowness) in $\Omega_j$ for $j=1,\ldots,K$. Here, $\mathbf{r}\in\Omega$ denotes the position coordinate and $\mathbf{p}$ denotes the momentum (or slowness) vector. We assume that each $\Omega_j$, $j=1,\ldots,K$ is a convex polygon containing a homogeneous medium, that is, $\eta(\mathbf{r})=\eta_j$ when $\mathbf{r}\in\Omega_j$ for $j=1,\ldots,K$ and $\eta_j$, $j=1,\ldots,K$ are constants. Note that when $\eta_j=c_j^{-1}$, the inverse of the phase velocity in $\Omega_j$, then
$H_j$ corresponds to the Hamiltonian for the ray trajectories
obtained in the GO limit for the Helmholtz equation
\begin{equation}\label{eq:Helm}
\Delta u + k_j^2 u =0,
\end{equation}
with $k_j=\omega/c_j$ the wavenumber at angular frequency
$\omega$. Alternatively when $\eta_j=\sqrt{\omega}/c_j$, then $H_j$ corresponds to ray trajectories for the frequency domain biharmonic wave equation
\begin{equation}\label{eq:Bih}
\Delta^2 u - k_j^4 u =0,
\end{equation}
in the GO limit. The additional $\sqrt{\omega}$ factor stems from the fact that the phase velocity $c_j\propto\sqrt{\omega}$ for (flexural) wave solutions to (\ref{eq:Bih}), whereas $c_j$ is independent of $\omega$ for solutions to the Helmholtz equation - see Appendix \ref{sec:App} for further details.

\begin{figure}
\centering
\includegraphics[width=0.6\textwidth]{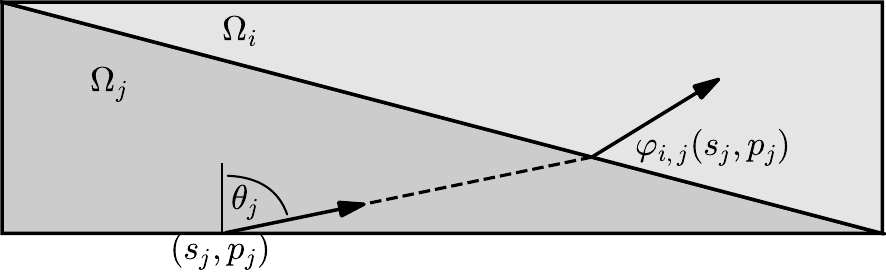}
\caption{The boundary map $\varphi_{i,j}(s_j,p_j)$ taking the
phase-space coordinates $X_j=(s_j,p_j)\in\Gamma_j\times(-\eta_j,\eta_j)$ to $\varphi_{i,j}(s_j,p_j)=(s'_i,p'_i)\in\Gamma_i\times(-\eta_i,\eta_i)$, which corresponds to the next intesection with a boundary edge where the ray undergoes either (refracted) transmission or a specular reflection.}\label{fig:BMap}
\end{figure}

We next define the phase-space coordinates $Y_j=(s_j,p_j)$ on the boundary $\Gamma_j=\partial\Omega_j$, $j=1,\ldots,K$,  where $s_j$ parametrizes $\Gamma_j$ by arclength and $p_j=\eta_j\sin(\theta_j)$ denotes the component of the momentum vector tangential to $\Gamma_j$ at $s_j$. The angle $\theta_j$ is formed between the outgoing ray trajectory and the inward pointing normal vector to $\Gamma_j$ at $s_j$ - see Fig.~\ref{fig:BMap}. Now let $\varphi_{i,j}(s_j,p_j)=(s'_i(s_j,p_j),p'_i(s_j,p_j))$ denote the boundary map describing the flow from the boundary of the domain $\Omega_j$ to $(s'_i(s_j,p_j),p'_i(s_j,p_j))$, where $s'_i(s_j,p_j)$ is positioned on the boundary of the domain $\Omega_i$ and $p'_i(s_j,p_j)=\eta_i\sin(\theta'_i(s_j,p_j))$ is the corresponding tangential momentum. Note that we are implicitly assuming that either $i=j$ or $\Gamma_i$ and $\Gamma_j$ share a common edge, through which the ray could travel. As before, the angle $\theta'_i$ is formed between the outgoing ray trajectory and the inward pointing normal vector to $\Gamma_i$ at $s'_i$. The calculation of the angle $\theta'_i$ will depend on the nature of the media within $\Omega_i$ and $\Omega_j$, as well as other physical factors such as the curvature in the case of a thin shell domain, for example. In the simplest possible cases, $\theta'_i$ will correspond to either a specular reflection if $i=j$ or the continued transmission of the ray in the same direction as when it arrives at $\Gamma_i$ if $i\neq j$, but $\eta_i=\eta_j$.

Phase-space densities are transported throughout $\Omega$ using a modified form of the FP operator (\ref{eq:FPO}) first proposed in \cite{GT09}, whereby the continuous time flow map $\varphi^\tau$ is replaced by the discrete boundary map $\varphi_{i,j}$. In this way the FP-operator becomes a local boundary integral operator $\mathcal{B}_{j}$, transporting a density $f$ from the phase-space on the boundary $\Gamma_j$ to the next boundary intersection with $\Gamma_i$ via
\begin{equation}\label{eq:OpB}
\mathcal{B}_{j}[f](X_i):= \int e^{-\mu_j D(X_i,Y_j)}w_{i,j}(Y_j) \delta(X_i-\varphi_{i,j}(Y_j)) f(Y_j)
\, \mathrm{d}Y_j.
\end{equation}
Here $X_i\in \Gamma_i\times(-\eta_i,\eta_i)$ for some $i=1,2,\ldots,K$ and the weighting factor $w_{i,j}$ is introduced to incorporate absorption factors at boundaries as well as direction-dependent reflection/transmission coefficients. A trajectory length dependent damping factor with coefficient $\mu_j\geqslant 0$ has also been introduced, with $D(X_i,Y_j)$ representing the Euclidean distance between $s_j$ and the solution point. The global boundary integral operator $\mathcal{B}=\sum_j\mathcal{B}_j$ is formed by taking the sum over each sub-domain $\Omega_j$ that shares a common edge with $\Omega_i$.

The stationary boundary density $\rho$  induced by a prescribed initial boundary density $\rho_{0}$ can be obtained from a Neumann series via
\begin{equation}\label{eq:RhoFinal}
\rho =
\sum_{n=0}^{\infty}\mathcal{B}^{n}[\rho_{0}] =
(I-\mathcal{B})^{-1}[\rho_{0}],
\end{equation}
where $\mathcal{B}^{n}$ models the density after $n$
iterates of the operator $\mathcal{B}$. Once the stationary density
$\rho$ is found, the interior density $\rho_{\Omega}$ can
then be obtained by projecting onto position space in $\Omega$
using
\begin{equation}\label{eq:RhoOmega}
\rho_{\Omega}(\mathbf{r}) = \frac{\eta_j^2}{\alpha} \int_{0}^{2\pi}
e^{-\mu_j D(\mathbf{r},s_j)}\rho(s_j(\mathbf{r},\Theta),p_j(\mathbf{r},\Theta))\,\mathrm{d}\Theta.
\end{equation}
The method for deriving (\ref{eq:RhoOmega}) from the underlying Hamiltonian mechanics is presented in Appendices B and C of Ref.~\cite{TH19}. The pre-factor $\alpha$ depends on the underlying wave equation, with $\alpha=1$ for the Helmholtz equation or $\alpha=2$ for the biharmonic equation. In addition, $\mathbf{r}\in\Omega_j$ is a prescribed solution point and
$\Theta\in[0,2\pi)$ is the polar angle parametrising trajectories approaching $\mathbf{r}$ from $s_j(\mathbf{r},\Theta)\in \Gamma_j$. With abuse of notation, the distance function $D$ is here used to represent the Euclidean distance between the solution point $\mathbf{r}\in \Omega_j$ and the boundary position $s_j\in\Gamma$.

\section{Discretization}\label{sec:FiniteDimAppr}
In this section we discuss a direction preserving discretization of the local boundary
operator (\ref{eq:OpB}) using a Petrov-Galerkin projection in order to numerically solve for the stationary density $\rho$ via (\ref{eq:RhoFinal}). We
detail how the two-dimensional integral can
be simplified and evaluated analytically for the case of convex polygonal sub-domains $\Omega_j$, $j=1,2,\ldots,K$ considered here.

\begin{figure}
\centering
\includegraphics[width=0.4\textwidth]{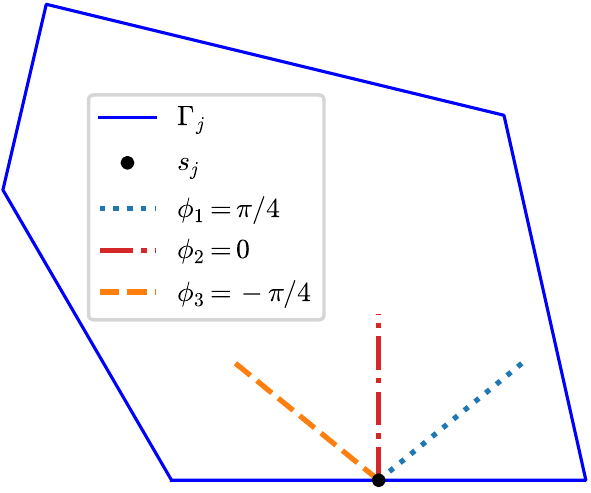}
\caption{Depiction of the local direction coordinates $\phi_n(s_j)\in(-\pi/2,\pi/2)$,  $n=1,2,3$ %
  for the case when there are $L=8$ global directions given by
  $\Phi_l=2\pi (l-1)/8$ for $l=1,2,\ldots, 8$. The map between local and global coordinates takes the form $\phi_n = \gamma - \Phi_{n+1}$ for $n=1,2,3$, where the global coordinate index has been shifted to run over the correct subset of directions. The change of sign between $\phi_n$ and $\Phi_{n+1}$ is due to the fact that the local angles are chosen positive for positive projected momenta $\tilde{p}_n$. The constant $\gamma$ corresponds to the global direction of the inward normal vector at $s_j$ and so $\gamma=\pi/2$ in the depicted example. 
}\label{fig:LocDir}
\end{figure}

Consider a subdivision of the boundary $\Gamma_j$ into elements $E^j_m$ for $m=1,2,\ldots,M_j$ and a set of global ray directions $\Phi_l\in[0,2\pi)$, where $l=1,2,\ldots L$, defined anti-clockwise relative to the positive $x$-axis. Here we make the choice $\Phi_l=2\pi (l-1)/L$, but note that this can be tailored for the example under
consideration to ensure the inclusion of dominant directional transmission paths if these are known \emph{a-priori}. We also perform the boundary element subdivision in such a way that any two sub-domains sharing a common edge will have identical boundary elements along the common edge. %
Furthermore, we set up the subdivision such that none of the $E^j_m$  extend over any of the vertices of the polygon $\Gamma_j$, i.e. every $E^j_m$ belongs to exactly one straight boundary edge.
Now define $\phi_n(s_j)\in(-\pi/2,\pi/2)$, $n=1,2,\ldots,N_m$ to be the local ray directions at $s_j\in\Gamma_j$. The local directions correspond to the subset of the global directions that are directed into $\Omega_j$ at $s_j$ and have been re-labelled according to the angle they make with the interior normal vector at $s_j$ - see Fig.~\ref{fig:LocDir} for a simple example case. %
If adjacent mesh cells are non-planar, then care has to be taken on
how to map the local to the global position set, see
Section~\ref{sec:flex-vibr-thin} for details.

We now approximate $\rho$ on $\Gamma_j\times(-\eta_j,\eta_j)$ using a finite-dimensional approximation of the form
\begin{equation}\label{eq:Expansion}
\rho(s_j,p_j)\approx\sum_{m=1}^{M_j}\sum_{n=1}^{N_m}\rho_{(j,m,n)}
b_{m}(s_j) \delta(p_j-\tilde{p}_n(s_j)), \qquad j=1,\dots,K,
\end{equation}
where $\tilde{p}_n(s_j)=\eta_j\sin(\phi_n(s_j))$ and $b_m(s_j)=|E^j_m|^{-1/2}$ for $s_j\in E^j_m$, and zero elsewhere, with $|E^j_m|=\text{diam}(E^j_m)$. For $m=1,2,\ldots, M_j$,  $b_m$ defines an orthonormal basis of piecewise constant functions with respect to the standard $L^2$ inner product and hence we apply a standard (Bubnov) Galerkin projection onto our basis in the position variable $s_j$. However, in the momentum variable $p_j$ we choose a set of test functions that are orthogonal (in fact orthonormal) in the $L^2$ inner product to $\delta(p_j-\tilde{p}_n(s_j))$ for $n=1,2,\ldots,N_m$ as follows. First define $\Delta\phi_n=\phi_{n+1}-\phi_n$ for $n=1,2,\ldots,N_m-1$ and consider a subdivision of the local direction range $$(-\pi/2,\pi/2)=\bigcup_{n=1}^{N_m}I_n$$ with $I_1=(-\pi/2,\phi_1+\Delta\phi_1/2]$,  $I_{N_m}=(\phi_{N_{m}-1}+\Delta\phi_{N_{m}-1}/2,\pi/2)$ and $$I_n=(\phi_{n-1}+\Delta\phi_{n-1}/2,\phi_n+\Delta\phi_n/2]$$ for $n=2,3,\ldots,N_m-1$. \DC{We then define the test functions as $$\chi_n(p_j)=\tilde\chi_n(\arcsin(p_j / \eta_{j})),$$ where $\tilde\chi_n(\theta_j) $ are the characteristic functions $\tilde\chi_n(\theta_j) = 1$ if $\theta_j\in I_n$, or zero otherwise.}
 Hence we have the property
\begin{equation}\label{eq:PetOrth}
\langle\delta(\cdot -\tilde{p}_n(s_j)),\chi_{n'}\rangle_{L^2(-\eta_j,\eta_j)}=0
\end{equation}
unless $n=n'$, when the inner product is one.

\begin{figure}
\centering
\includegraphics[width=0.5\textwidth]{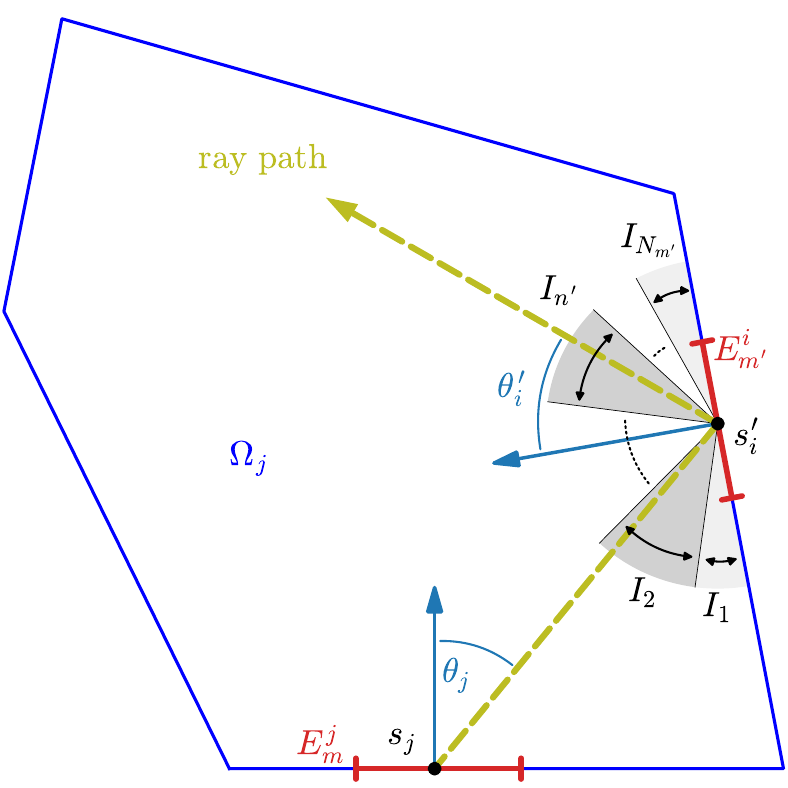}
\caption{Illustration of the case when the entries of the discretized boundary operator matrix $B_{I,J}$ are non-zero for a reflected ray ($i=j$). The initial position $s_j$ and local direction $\theta_j$ must be such that the ray arrives at the boundary element $E^i_{m'}$ and the local direction of the reflected (or transmitted if $i\neq j$) ray must fall within the local direction sub-interval $I_{n'}$.}\label{fig:Bij}
\end{figure}

A Petrov-Galerkin projection of the operator $\mathcal{B}_{j}$ on to the basis and test function combination described above leads to a matrix representation $B$ with entries given by
\begin{align}\label{eq:TransferBij}
\begin{split}
B_{I,J}&=\int_{\Gamma_j\times(-\eta_j,\eta_j)}\hspace{-13mm}e^{-\mu_j D(\varphi_{i,j}(Y_j),Y_j)} w_{i,j}(Y_j)b_{m}(s_j)b_{m'}(s'_i(Y_j))
\delta(p_j-\tilde{p}_n(s_j))\chi_{n'}(p'_i(Y_j))\,\mathrm{d}Y_j\vspace{2mm}\\
&=\int_{\Gamma_j}\hspace{-1mm} e^{-\mu_j D_i(s_j)} w_{i,j}(s_j,\tilde{p}_n(s_j))b_{m}(s_j)b_{m'}(s'_i(s_{j},\tilde{p}_n(s_j)))
\chi_{n'}(p'_i(s_{j},\tilde{p}_n(s_j)))\,\mathrm{d}s_j\\
&=\frac{1}{|E^j_m|^{1/2}}\int_{E^j_m}\hspace{-1mm}e^{-\mu_j D_i(s_j)} b_{m'}(s'_i(s_{j},\tilde{p}_n(s_j)))
\chi_{n'}(p'_i(s_{j},\tilde{p}_n(s_j)))\,\mathrm{d}s_j,\end{split}
\end{align}
where $I$ and $J$ denote the multi-indices $I=(i,m',n')$ and
$J=(j,m,n)$, respectively. We have also introduced the notation $D_i(s_j)$ for the Euclidean distance between $s_j\in\Gamma_j$ and $s_i'(s_j,\tilde{p}_n(s_j))\in\Gamma_i$. Note that the four dimensional integration in the definition of the Galerkin projection has been reduced to a single integral over the boundary element $E^j_m\subset\Gamma_j$ due to the local support of the spatial basis and the properties of the Dirac $\delta$ distributions arising in equations (\ref{eq:OpB}) and (\ref{eq:Expansion}). To obtain the third line in (\ref{eq:TransferBij}), we have further assumed that the weight function $w_{i,j}$ is independent of $s_j\in E^j_m$ and may define either locally constant damping at boundaries or direction-only dependent reflection/transmission coefficients.

The calculation of the matrix element $B_{I,J}$ is therefore relatively simple, since the two basis functions are also locally constant and will be zero unless the direction $\theta'_i\in I_{n'}$ and $s_i'\in E^i_{m'}$, meaning that most entries to the matrix will be zero. The case when the matrix elements are non-zero is illustrated in Fig.~\ref{fig:Bij} and the calculation only involves the integral of the exponential term $e^{-\mu_j D_i(s_j)}$ over the element $E^j_m$ and multiplication by the pre-factor $ w_{i,j}(\tilde{p}_n(s_j))b_m(s_j)b_{m'}(s'_i)=w_{i,j}(\tilde{p}_n(s_j))(|E^j_m||E^i_{m'}|)^{-1/2}$. For polygonal boundaries, the Euclidean distance function $D_i(s_j)$ is linear in $s_j\in E^j_m$ and hence the integral in the third line of (\ref{eq:TransferBij}) can be performed analytically with relative ease. The only potential complication arises when $s_i'(s_j,\tilde{p}_n(s_j))$ coincides with one of the vertices, and the integral must then be sub-divided at the corresponding value of $s_j$. We note that analytic spatial integration is also possible for higher order basis functions using Legendre polynomials, albeit with a more complicated implementation process based on recursion formulae \cite{JB17}. Since the expansion (\ref{eq:Expansion}) represents a local density approximation and the transfer of energy is only between connected sub-domains, $B$ is also a block-sparse matrix for large $K$.

The coefficients of the
expansion (\ref{eq:Expansion}) can be found by solving the linear
system ${\boldsymbol{\rho}} = (I-B)^{-1}{\boldsymbol{\rho}}_0$,
which corresponds to the discretized form of equation
(\ref{eq:RhoFinal}). Here $\boldsymbol{\rho}_0$ and
$\boldsymbol{\rho}$ represent the coefficients of the expansions of
$\rho_{0}$ and $\rho$, respectively, when
projected onto the finite dimensional basis. Note that for $I-B$ to
be invertible, we must include either a positive dissipation factor $\mu_j>0$ or a boundary absorption factor within $w_{i,j}$. The entries of the source vector
$\boldsymbol{\rho}_0$ corresponding to an initial density $\rho_0$ are given using the property (\ref{eq:PetOrth}) and orthonormality of the spatial basis functions via
\begin{align}\label{eq:SourceDisc}
\begin{split}
[\boldsymbol{\rho}_0]_J&=\int_{\Gamma_j\times(-\eta_j,\eta_j)}\rho_0(s_j,p_j)b_m(s_j)\chi_{n}(p_j)\mathrm{d}Y_j\\
&=\frac{1}{|E^j_m|^{1/2}}\int_{E^j_m}\int_{I_n}\rho_0(s_j,p_j)\mathrm{d}p_j\mathrm{d}s_j.
\end{split}
\end{align}

Once $\boldsymbol{\rho}$ has been computed and
substituted into (\ref{eq:Expansion}), then the interior density
$\rho_\Omega$ can be approximated using (\ref{eq:RhoOmega}) as follows
\begin{align}\label{eq:RhoOmegaD}
\begin{split}
\rho_{\Omega}(\mathbf{r}) &\approx \\
&\frac{\eta_j^2}{\alpha} \sum_{m=1}^{M_j}\sum_{n=1}^{N_m}\rho_{(j,m,n)}\int_{0}^{2\pi}\hspace{-2mm}e^{-\mu_j D(\mathbf{r},s_j(\mathbf{r},\Theta))}
\delta(p_j(\mathbf{r},\Theta)-\tilde{p}_n(s_j(\mathbf{r},\Theta)))b_m(s_j(\mathbf{r},\Theta))\,\mathrm{d}\Theta.
\end{split}
\end{align}
Applying the properties of the Dirac delta one may write
\begin{equation}\label{eq:IntDelta}
  \delta\left(p_j(\mathbf{r},\Theta)-\tilde{p}_n(s_j(\mathbf{r},\Theta)) %
    \rule{0pt}{2.1ex}\right)=\frac{\delta(\theta_j(\mathbf{r},\Theta)-\phi_n)}{\eta_j\cos(\theta_j(\mathbf{r},\Theta))}=\frac{\delta(\Theta-\Phi_l)}{\eta_j\cos(\theta_j(\mathbf{r},\Theta))},
\end{equation}
where $\Phi_{l}$ is the global direction corresponding to the local direction $\phi_n$,  which arises in the term $\tilde{p}_n(s_j(\mathbf{r},\Theta))=\eta_j\sin(\phi_n(s_j(\mathbf{r},\Theta)))$. We note that the $s_j$ dependence of $\phi_n$ relates to the differences in the local direction set on each edge of $\Gamma_j$, and is independent of $s_j$ along a given edge of $\Gamma_j$. To apply the  expression (\ref{eq:IntDelta}) we must therefore subdivide the integral in (\ref{eq:RhoOmegaD}) into a set of sub-integrals, split at the angles $\Theta$ where $s_j(\mathbf{r},\Theta)$ corresponds to a vertex of $\Gamma_j$. The second equality in (\ref{eq:IntDelta}) is due the fact that the mapping between local and global directions is a switch of orientation (multiplication by $-1$) and then a translation by a local edge dependent constant - see Fig.~\ref{fig:LocDir} and caption. The denominator in (\ref{eq:IntDelta}) is calculated from the (absolute value of) the $\theta_j$ derivative of $p_j(\mathbf{r},\Theta)=\eta_j\sin(\theta_j(\mathbf{r},\Theta))$. Finally, we arrive at the following global ray summation result for the interior density
\begin{equation}\label{eq:RhoOmegaDS}
\rho_{\Omega}(\mathbf{r}) \approx \frac{\eta_j}{\alpha} \sum_{l=1}^{L}
\frac{e^{-\mu_j D(\mathbf{r},s_j(\mathbf{r},\Phi_l))}\rho_{(j,m,n)}(\mathbf{r},\Phi_l)}{|E_m^j(\mathbf{r},\Phi_l)|^{1/2}\cos(\theta_j(\mathbf{r},\Phi_l))}.
\end{equation}
Here the integral has only returned a non-zero value each time $\Phi$ coincides with a member of the global direction set and so the double sum over boundary elements and local directions can be replaced with a single summation over the global directions $\Phi_l$, $l=1,2,\ldots, L$. The boundary position $s_j$ and local direction $\theta_j$ can be determined from knowledge of the solution point $\mathbf{r}$ and the ray direction $\Phi_l$, hence we can also obtain the element index $m$ and the local direction index $n$.
\section{Numerical results}\label{sec:NumResults}
In this section we study a series of examples in order to investigate the effectiveness of the proposed discretization, including rates of convergence as the phase-space discretization is refined and investigating the effect of including a transparent internal mesh. We include examples both where the directivity can be exactly described by a finite direction set and irregular domains where the directivity can only be approximated by increasing the number of ray directions. We also consider both flat two-dimensional acoustic domains and flexural vibrations of curved thin shells. We start with the simplest example of a ray reflecting back and forth between two parallel walls.

\subsection{Simple parallel wall reflections}\label{sec:simpeg}

Consider a unit square acoustic domain $\Omega$ with vertices $(0,0)$, $(1,0)$, $(1,1)$ and $(0,1)$ and sound-hard reflections along the adjoining edges. The problem is forced via a spatially constant line source
\begin{equation}\label{eq:SqRho0}
\rho_0(s,p)=\mathbf{1}_{\{x_s=0\}}(s)\delta(p)
\end{equation}
along the fourth edge from $(0,1)$ to $(0,0)$, propagating rays into the domain directed perpendicularly to this edge. Note that since there is only one sub-domain we have dropped the subscripts from the phase-space boundary coordinates and denote $(s,p)=(s_1,p_1)$. The notation $\mathbf{1}_{\{x_s=0\}}$ is used to represent an indicator function that is zero unless the  position defined by the arclength parameter $s$ has (Cartesian) $x$-coordinate $x_s=0$, where $\mathbf{1}_{\{x_s=0\}}(s)=1$.

The resulting ray dynamics is extremely simple since rays only bounce between the source edge and the opposite edge, and are always directed perpendicular to these edges. The ray density will simply take a constant value on these two edges and be zero on the remaining two edges. After applying a dissipative term of the form $\exp(-\mu D)$, with damping coefficient $\mu$ and trajectory length $D$, then the interior ray density $\rho_{\Omega}$ may be calculated at $\mathbf{r}=(x,y)\in\Omega$ via a geometric series over the reflection order to give
\begin{equation}\label{eq:SqEx}
\rho_\Omega (x,y) = \frac{e^{-\mu x}+e^{-\mu(2-x)}}{1-e^{-2\mu}}.
\end{equation}
For the calculations in this section we set the damping parameter $\mu=\pi/2$.

Despite its apparent simplicity, we study this example to demonstrate the proposed methodology in comparison to previous DEA approaches which employ the same spatial discretization as here, but instead use a Bubnov-Galerkin projection onto a basis of scaled Legendre polynomials in the momentum variable \cite{CTLS13, TH19}. We first consider this problem as a single domain and then as a multi-domain problem, where the sub-domains are given by triangular mesh cells. We then introduce curvature in the $x$-direction and consider essentially the same example with the curvilinear coordinate taking the role of $x$ in the exact solution.

\begin{table}
  \caption{\noindent Relative mean errors and estimated orders of convergence as the size of the momentum basis is varied: comparison between the best and worst choice of global direction set for the Dirac delta basis and an orthogonal polynomial basis using Legendre polynomials.}
    \label{tab:SqDelt}
  \centering
  \begin{tabular}{|c|c|c|c|c|c|} \cline{2-6}
   \multicolumn{1}{c}{} & \multicolumn{1}{|c|}{$\Phi_l=2\pi (l-1)/L$} & \multicolumn{2}{|c|}{$\Phi_l= (2l-1)\pi/L$} & \multicolumn{2}{|c|}{Legendre basis}\\ \hline
  $L/2$ or basis order & Error & Error & EOC & Error  & EOC \\ \hline
4 & 1.4737e-16  & 6.2325e-2 & -     & 5.6847e-2 & -           \\
8 & 1.4736e-16  & 5.8211e-2 & 0.10  & 3.4446e-2 & 0.72        \\
16 & 1.4736e-16 & 1.4781e-2 & 1.98  & 2.0904e-2 & 0.72        \\
32 & 1.4736e-16 & 7.7382e-3 & 0.93  & 1.1887e-2 & 0.81        \\
64 & 1.4736e-16 & 3.8344e-3 & 1.01  & 6.4164e-3 & 0.89        \\
128 & 1.4736e-16 & 1.4887e-3 & 1.37  & 3.5660e-3 & 0.85        \\
256 & 1.4736e-16 & 3.5971e-4 & 2.05  & 1.8413e-3 & 0.95         \\   \hline
  \end{tabular}
\end{table}

Considering the unit square as a single domain, we divide the boundary into four elements (one per edge) and apply a direction basis containing four directions $\Phi_1=0$, $\Phi_2=\pi/2$, $\Phi_3=\pi$  and $\Phi_4=3\pi/2$. Using this coarse discretization we are able to achieve a machine accuracy approximation of $\rho_{\Omega}$. The reason for this effectiveness is that even at this coarse discretization level, the exact solution for $\rho$ on $\Gamma$ lies in the approximation space since it takes a constant value on each edge and propagates only in the directions $\Phi_1$ and $\Phi_3$. One can also consider the convergence of the model if the direction set is chosen to be an equi-spaced set of directions that are as far from the true propagation directions as possible. That is, we replace the choice $\Phi_l=2\pi (l-1)/L$  with $\Phi_l= (2l-1)\pi/L$, $l=1,2,\ldots,L$. The results are summarised in Table \ref{tab:SqDelt} and also include a comparison against Legendre basis approximation typically used in DEA simulations. The first column relates to the number of degrees of freedom in the momentum approximation on a given boundary element, which for the Dirac delta basis is $L/2$ since half of the equally spaced global direction set are designated as directions propagating within the unit square domain. For the Legendre basis, the number of degrees of freedom is actually the basis polynomial order (shown in the first column) plus one since the basis order starts from zero. The error columns show the relative mean error sampled at $P=3638$ interior points $\mathbf{r}_i$, $i=1,2,\ldots,P$ taken as the centroids of a triangle mesh generated by the Distmesh package \cite{distmesh} for MATLAB with mesh spacing 0.025 (corresponding to the finest internal mesh considered in the next experiment). Explicitly, we calculate
\begin{equation}\label{MRE}\text{Error}=\frac{\sum_{i=1}^{P}|\hat{\rho}_\Omega(\mathbf{r}_i)-\rho_\Omega(\mathbf{r}_i)|}{\sum_{i=1}^{P}\rho_\Omega(\mathbf{r}_i)},
\end{equation}
where $\hat{\rho}_\Omega$ is the approximation to $\rho_\Omega$ calculated via (\ref{eq:RhoOmegaDS}) for the error values reported in columns 2 and 3, or calculated using the DEA approximation detailed in Refs.~\cite{CTLS13, TH19} for the error values in column 5.
The estimated order of convergence (EOC) is then calculated via the binary logarithm of the ratio of each error value to the previous error value (in the same column).

The results in Table \ref{tab:SqDelt} show that even in the worst case scenario, the Dirac delta basis achieves reasonably small errors but with a fairly volatile convergence rate that is approximately between first and second order. For any Dirac delta basis including the true propagation directions, the method is accurate to machine precision. In comparison to a DEA approximation using Legendre polynomials, even the worst case direction choice typically provides better accuracy when compared to the Legendre polynomial basis with one additional degree of freedom. Whilst the convergence rate for the Legendre polynomial basis is less volatile, it is also typically sub-linear and slower than the convergence of the Dirac delta basis with the `bad' direction set choice. Two further advantages of the Dirac delta basis are that one achieves a sparser representation of the matrix $B$ and a simpler evaluation of the matrix entries, which is free of the numerical evaluation of the integral with respect to the momentum variable due to the sifting property of the Dirac delta. \DC{The computation times approach an $\mathcal{O}(L)$ scaling from below as the number of directions is increased. For non-optimised serial MATLAB code running on a 2.6GHz processor, the computational times for the Dirac delta basis results in Table \ref{tab:SqDelt} range between 8s and 255s, and in each case the time is dominated by the assembly cost for the linear system. All other tasks (preprocessing, solution of the linear system, calculation of the interior density) take around 2s in total, for all values of $L$ considered. We note that this approximate $\mathcal{O}(L)$ scaling is a consequence of the sparsity of the matrix (\ref{eq:TransferBij}). The number of matrix entries to compute for the linear system formed using the Legendre basis grows with the square of the basis order and the computational cost scales worse than this because higher order basis functions require additional numerical integration effort.}  

\begin{figure}
\centering
\includegraphics[trim=60 75 20 75, width=0.85\textwidth]{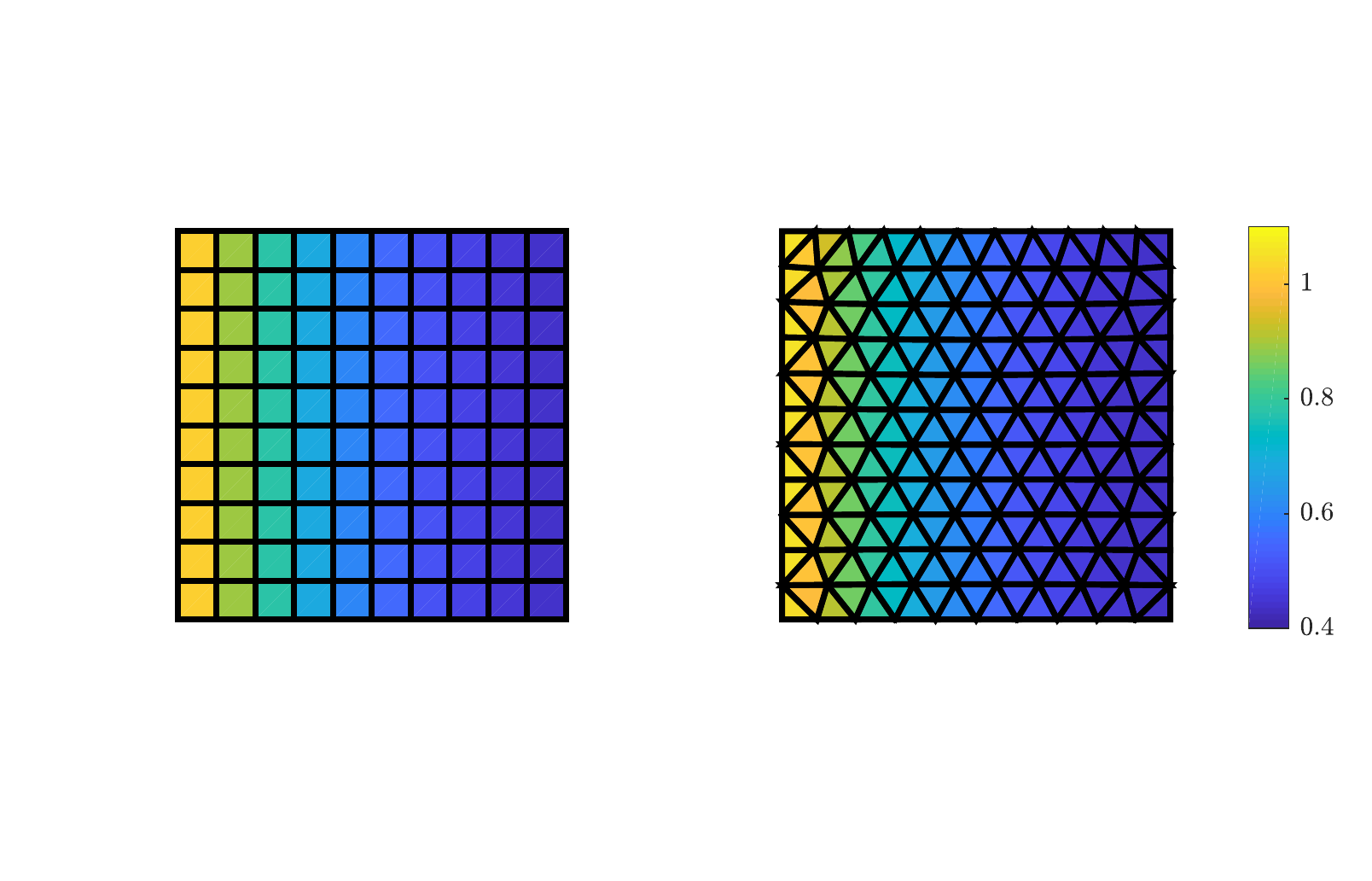}
\caption{The exact solution (\ref{eq:SqEx}) illustrated on the two types of internal mesh considered. Left: mesh consisting of $K=100$ square mesh cells generated from an equal subdivision of the standard Cartesian coordinates, right: mesh consisting of $K=211$ triangles generated by Distmesh with nominal mesh spacing $h=0.1$.}\label{fig:Meshes}
\end{figure}

We now consider the effect of including an internal mesh in the model and modelling the ray flow as a local propagation between mesh cell boundaries. This procedure is in no way necessary for the example at hand, but would be necessary to describe propagation through complex domains including curvature and/or inhomogeneities, and so will be tested here. We consider two types of internal mesh as depicted in Fig.~\ref{fig:Meshes}, the first of which is based on square mesh cells resulting from an equal subdivision in both $x$ and $y$. The second is a triangulation generated by the Distmesh package for MATLAB \cite{distmesh}. In each case we apply the simplest possible boundary element discretization on each mesh cell $\Omega_j$, $j=1,2,\ldots, K$ with exactly one boundary element per edge. To consider the effect of the mesh on the accuracy in isolation of other factors we choose the direction set $\Phi_l= 2\pi(l-1)/L$ with $L=4$, which calculated the result to machine accuracy in the absence of an internal mesh.

\begin{table}
  \caption{\noindent Relative mean errors and estimated orders of convergence as the size of the internal mesh elements are varied: comparison between square and triangular elements for the Dirac delta basis and triangular elements for different orders of Legendre polynomial basis.}
    \label{tab:IntMesh}
  \centering
  \begin{tabular}{|c|c|c|c|c|c|c|c|} \cline{2-8}
   \multicolumn{1}{c}{} & \multicolumn{1}{|c|}{Quad delta} & \multicolumn{2}{|c|}{Distmesh delta} & \multicolumn{2}{|c|}{Distmesh poly 8}  & \multicolumn{2}{|c|}{Distmesh poly 256}\\ \hline
  $h$ & Error  & Error & EOC & Error  & EOC & Error  & EOC \\ \hline
0.4 & -              & 3.9038e-2 & -     & 6.2347e-2 & -   & 4.8582e-2 & -           \\
0.2 & 3.0060e-16    & 1.3012e-2 & 1.59  & 5.4538e-2 & 0.19 & 1.6165e-2 & 1.59        \\
0.1 & 1.4942e-15    & 5.8543e-3 & 1.15  & 1.0237e-1 & -0.91 & 1.0011e-2 & 0.69        \\
0.05 & 5.9241e-15  & 2.7935e-3 & 1.07  & 1.2036e-1 & -0.23 & 5.6531e-3 & 0.82        \\
0.025 & 2.5038e-14  & 1.3721e-3 & 1.03  & 1.2070e-1 & 0.00 & 4.5696e-3 & 0.31        \\ \hline
  \end{tabular}
\end{table}

The first column of Table \ref{tab:IntMesh} shows the mesh spacing, which is either the length of each side of the square element mesh or the nominal mesh spacing parameter used within the Distmesh package. The remaining columns show the corresponding errors and EOCs calculated as in Table \ref{tab:SqDelt}, except now  taking the average error using the centroids of the mesh cells within the current internal mesh (as opposed to the centroids of the fixed Distmesh triangulation with $h=0.025$). The second to fourth columns show results computed using the Dirac delta basis with the global direction set $\Phi_l= 2\pi(l-1)/L$ with $L=4$ and the final four columns compare these results to those computed using a Legendre polynomial approximation truncated after the 8th or 256th order term as indicated in the table headings. The results in the 2nd column show that introducing the quadrilateral mesh retains the machine precision accuracy obtained in the absence of a mesh, but as the mesh density is increased, round off errors accumulate and become more significant. The reason for the high accuracy is that the exact solution (\ref{eq:SqEx}) is independent of the $y$-coordinate, meaning that it takes a constant value along the vertical boundary lines of each mesh cell. Furthermore, since the only rays propagating through $\Omega$ travel parallel to the $x$-axis, then rays do not emerge from the horizontal boundary lines of a mesh cell making the boundary density $\rho$ zero along these edges. Hence, the exact solution for $\rho$ on $\Gamma$ is again within the approximation space. This is no longer the case for general triangle meshes generated using Distmesh and consequently one loses the high accuracy as shown in column 3 of Table \ref{tab:IntMesh}. The error introduced corresponds to the approximation of the exact solution (\ref{eq:SqEx}) by piecewise constant functions along the mesh cell boundaries within $\Omega$ and consequently the EOC results shown in the fourth column of Table \ref{tab:IntMesh} indicate first order convergence. These results compare favourably with the Legendre basis results. Here, the refinement of the spatial mesh with a {\em fixed} order polynomial approximation in momentum leads to convergence to a solution different from the exact solution; this solution becomes closer to the exact solution only once the order polynomial approximation in the momentum variable is increased.

\begin{figure}
\includegraphics[trim=50 150 50 0,width=0.6\textwidth]{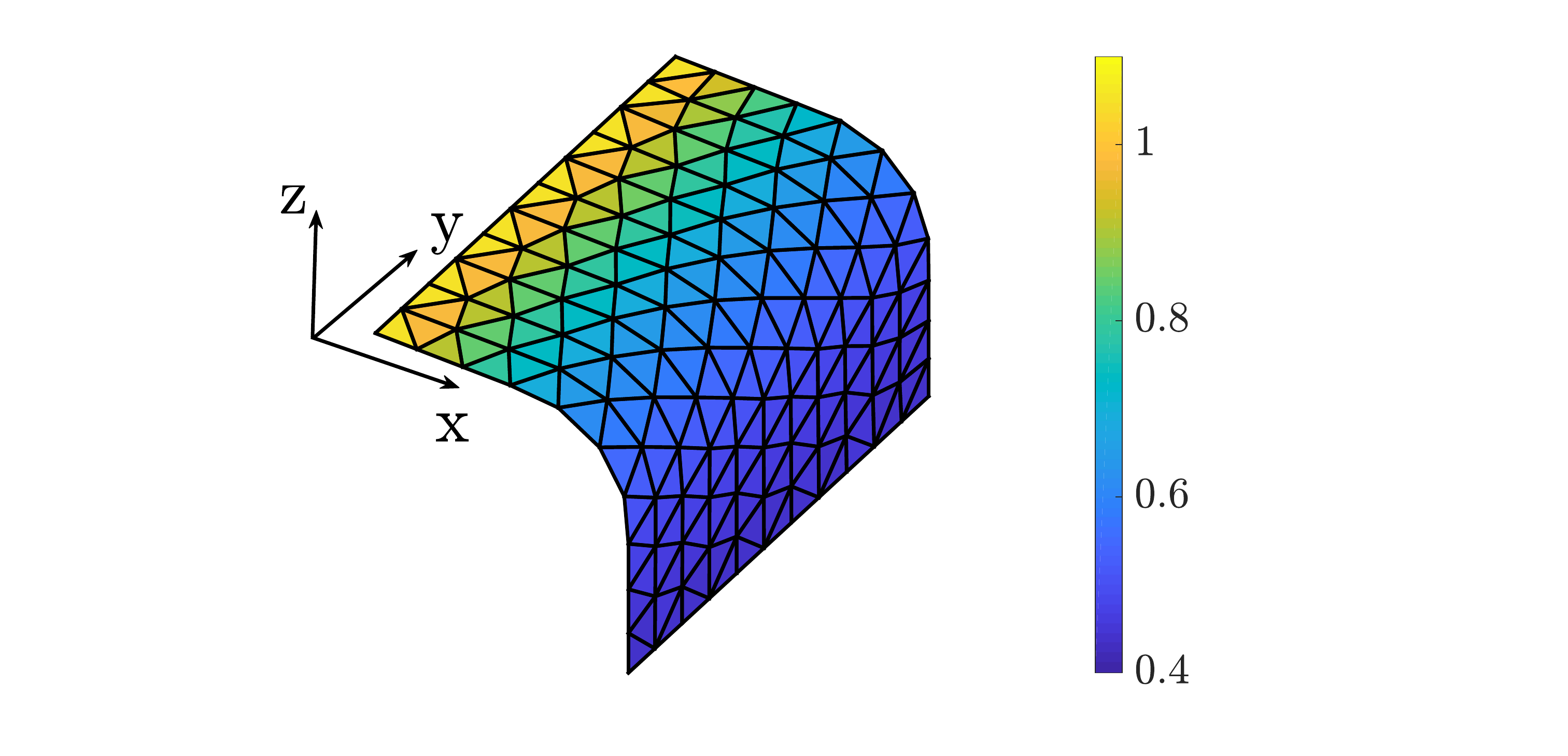}
\hfill\begin{tabular}{|c|c|c|} \cline{2-3}
   \multicolumn{1}{c}{}  & \multicolumn{2}{|c|}{Distmesh delta} \\ \hline
  $h$ & Error & EOC   \\ \hline
0.4          & 8.6884e-2 & -      \\
0.2    & 1.8460e-2 & 2.23         \\
0.1     & 6.8713e-3 & 1.43        \\
0.05   & 2.9971e-3 & 1.20         \\
0.025   & 1.4194e-3 & 1.08        \\ \hline
\end{tabular}
\caption{ %
  \label{fig:Ridge0} %
  \noindent Results for the parallel wall reflection problem when the middle third of the domain in the $x$-direction has been replaced with a quarter cylinder with the same surface area. Left: numerical solution computed using the Dirac delta basis containing the propagation directions of the exact solution and with one spatial (boundary) element per mesh cell edge using a mesh with $h=0.1$. Right: table of results showing the relative mean errors and estimated orders of convergence as the size of the internal mesh is varied.}
\end{figure}

We conclude this section by considering the additional error introduced when the triangulated domain is used to approximate a curved surface. The parallel wall reflection problem and boundary conditions are unchanged from before, except that the domain is now curved in the $x-z$ plane. The problem retains invariance in the $y$-direction. A curvilinear coordinate $\tilde{x}\in[0,1]$ takes the role of the $x$-coordinate in the exact solution (\ref{eq:SqEx}) and the definition of the global direction set. The region $\tilde{x}\in(1/3,2/3)$ is a quarter cylinder with radius $2/(3\pi)$ and outside this region the domain is flat as before - see Fig.~\ref{fig:Ridge0}. The additional error in this example (compared with the results shown in columns 3 and 4 of Table \ref{tab:IntMesh}) arises solely from the approximation of the curved geometry in the region $\tilde{x}\in(1/3,2/3)$ by a mesh of flat triangles. A ray travels a shorter distance through the triangle mesh than it does on the actual curved surface. The distance travelled through the mesh will converge towards the correct one as the mesh is refined.

 The numerical solution on a Distmesh generated mesh with $h=0.1$ is shown in the left plot of Fig.~\ref{fig:Ridge0} and we note the correspondence to the exact solution on the equivalent Distmesh generated mesh for the unit square domain shown in the right plot of Fig.~\ref{fig:Meshes}. The table at the right of the plot in Fig.~\ref{fig:Ridge0} shows the relative mean errors and EOCs when the mesh size is varied as in Table \ref{tab:IntMesh}. The results in columns 3 and 4 of Table \ref{tab:IntMesh} show the results for the corresponding examples on the unit square. Introducing the curved geometry has therefore resulted in larger errors than those for the equivalent flat planar domain problem as would be expected. However, the increase in the error is relatively small for $h\leq 0.1$ and the convergence rate appears to be unaffected, approaching first order as the mesh is refined. \DC{For the calculations in Table \ref{tab:IntMesh} and those documented in Fig.~\ref{fig:Ridge0}, the direction basis size $L$ was fixed and instead the number of mesh cells $K$ was varied. The computational cost scales with $\mathcal{O}(K^2)$ and for large $K$ it is recommended to compute in parallel, since the method is embarrassingly parallel in $K$ \cite{CTLS13}. As before, the computation times are dominated by the cost of assembling the linear system, but there is a notable increase in the preprocessing cost compared to before, owing to the mesh generation and associated pre-calculation tasks (finding nearest neighbours, classifying edges etc.). The computation times ranged from 0.3s to 34 minutes for non-optimised MATLAB code using the parallel computing toolbox on 16 cores with a 2.6GHz clock speed.}

\subsection{Irregular coupled cavity examples}\label{sec:CCnum}

\begin{figure}
\centering
\includegraphics[trim= 140  0 140  0, width=\textwidth]{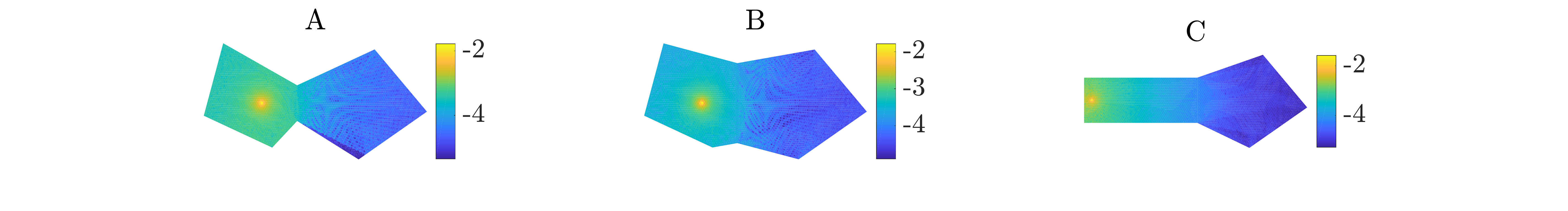}
\caption{Logarithm (base 10) of the interior energy densities in coupled cavity configurations A, B and C with damping coefficient $\mu=\pi/2$. The results have been calculated using $L=2048$ global ray directions and an average boundary element step size of $0.0125$.}\label{fig:ConABC}
\end{figure}

We now consider numerical examples in three different configurations (labelled A to C) of coupled polygonal acoustic cavities $\Omega=\Omega_1\cup\Omega_2$ as introduced in Refs.~\cite{BM05, GT09} and shown in Fig.\ \ref{fig:ConABC}. In each case the cavities are connected along the line $y=0$ and driven by an acoustic velocity potential point source with angular frequency $\omega$ located at $\mathbf{r}_0\in\Omega_1$. This leads to an initial boundary density in $\Omega_j$, $j=1,2$, of the form \cite{JB17}
\begin{equation}\label{eq:Rho0pt}
\rho_0(s_j,p_j;\mathbf{r}_0)=\frac{\omega\varrho\cos(\theta_0)e^{-\mu D(\mathbf{r}_0,s_j)}w_{j,1}(s_j,p_j)\delta(p_j-p_0)}{8\pi D(\mathbf{r}_0,s_j)},
\end{equation}
where $\varrho$ is the density of the acoustic fluid in $\Omega$. In addition, $p_0=\eta_1\sin(\theta_0)$ is the tangential momentum of the ray emerging from $s_j$ (after specular reflection in $\Omega_1$ or transmission into $\Omega_2$) that arrived directly from $\mathbf{r}_0$, and $\theta_0\in(-\pi/2,\pi/2)$ is the angle that this ray makes with the internal normal vector at $s_j$. The initial boundary density $\rho_0$ will only be non-zero along polygonal edges that are directly illuminated by rays from the source point, meaning that $\rho_0(s_2,p_2)=0$ along the edges of $\Gamma_2$ except for the edge that connects to $\Omega_1$. This property enters (\ref{eq:Rho0pt}) through the term
$$
w_{j,1}(s_j,p_j)=\left\{\begin{array}{ll}1 & \text{if }j=1\text{ and }s_j\text{ is on a free edge,}\\
R(p_j) & \text{if }j=1\text { and }s_j\text{ is on the interface connecting to }\Omega_2,\\
T(p_j) & \text{if }j=2\text{ and }s_j\text{ is on the interface connecting to }\Omega_1,\\
0 & \text{otherwise.}
\end{array}\right.
$$
Here, by a free edge we mean an edge not connecting to another sub-domain and $R$/$T$ denote the reflection/transmission probabilities for rays in $\Omega_1$ arriving at the interface between $\Omega_1$ and $\Omega_2$. In the examples here we assume that $\Omega$ contains a homogeneous fluid medium meaning that $R(p_j)\equiv 0$ and $T(p_j)\equiv 1$; furthermore, we choose: $\varrho=1$, $c_j=1$ for $j=1,2$, $\omega=100\pi$ and $\mu=\pi/2$, as before.

Figure \ref{fig:ConABC} shows the results for $\log_{10}(\rho_\Omega)$ in each of the coupled cavity configurations using $L=2048$ global ray directions and an average boundary element step size of $0.0125$, that is the number of elements along each edge is calculated by dividing the edge length by $0.0125$ and rounding to the nearest integer. The interior points at which the solution is calculated are taken as the centroids of a triangle mesh generated by the Distmesh package \cite{distmesh} for MATLAB with mesh spacing 0.025 as before.  For each configuration one observes that the interior solutions fluctuate locally between illuminated and darker regions, particularly in the receiver sub-domain $\Omega_2$. This can be attributed to the momentum discretization using a finite number of discrete ray directions and such artefacts are typical in ray tracing methods - see, for example, Section IV(B) of Ref.~\cite{LS15}. Note that the fluctuations are less visible in the source sub-domains since the solution in $\Omega_1$ is the sum of the interior density computed via finite ray summation (\ref{eq:RhoOmegaDS}) and the initial acoustic energy density emitted from the source before reaching the boundary, which is given by $\varrho\omega^2\eta_1^2|G(\mathbf{r},\mathbf{r}_0)|^2$, where $G$ is the free-space Green's function for the Helmholtz equation. Note that (\ref{eq:RhoOmegaDS}) only gives the solution after the source has arrived at the boundary meaning that the initial direct source illumination has to be added to give the full solution.

We now compare our results to an `exact' solution given by the image source method~\cite{JC09,LS15}, which is exact in the sense that we perform sufficiently many iterations (reflections) that the mean solution in each sub-domain has converged to 12 significant digits. The number of iterations required for this level of convergence depends on the damping coefficient $\mu$ and for the value $\mu=\pi/2$ used here, the solution converges up to 12 digits after 32 iterations.  Calculating the relative errors with respect to this exact solution using (\ref{MRE}) gives poor results due to the fluctuations described above. For the discretizations used to create Fig.~\ref{fig:ConABC}, the relative errors for the solutions in the receiver sub-domain $\Omega_2$ range from $0.083338$ for configuration C, to $0.20644$ in configuration B and up to $0.25379$ for configuration A. This is considerably worse than expected, given the results in the previous section. The more accurate results for configuration C are probably due to the source sub-domain $\Omega_1$ being rectangular, meaning that here the reflected ray directions will be exactly represented in the direction basis.

\begin{figure}
\centering
\includegraphics[trim= 40  0 40  0, width=\textwidth]{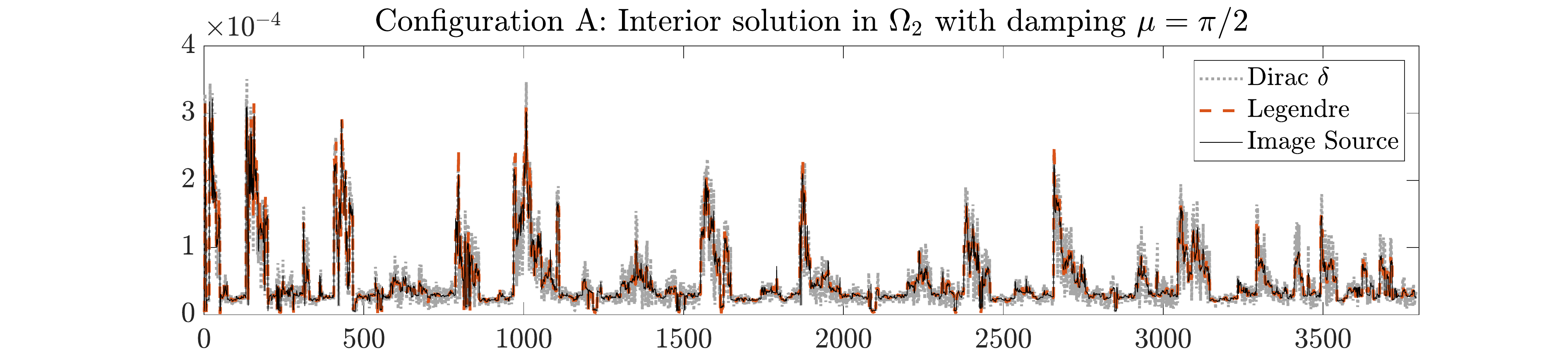}
\caption{Interior energy density in configuration A calculated at $P=3791$ points in the right-hand (receiver) sub-domain $\Omega_2$ with damping coefficient $\mu=\pi/2$. The plot compares the results calculated using a Dirac delta momentum basis with $L=2048$ global ray directions and an average boundary element step size of $ 0.0125$, DEA with a Legendre polynomial expansion up to degree 32 in the momentum variable and the same spatial approximation as before, and the image source method with 32 iterations.}\label{fig:ConA}
\end{figure}

To better understand the above results, in Fig.~\ref{fig:ConA} we plot the interior density in the receiver sub-domain for the worst case (configuration A) against the interior point index number and compare to the image source result. Figure \ref{fig:ConA} also shows a corresponding DEA result computed using the same spatial approximation with average boundary element step size of $0.0125$, but with the approximation in the momentum variable computed using a Legendre polynomial expansion up to degree 32. The Legendre polynomial basis leads to a less fluctuating solution and the relative error here is $0.067113$. While this result is a considerable improvement on the Dirac $\delta$ basis approximation and uses fewer degrees of freedom in the momentum variable, it is still relatively large compared to the results of the previous section. The Dirac $\delta$ basis result in Fig.~\ref{fig:ConA} appears to approximately fluctuate about a mean given by the image source method result. To verify whether this is indeed the case we calculate the relative error in the mean of the Dirac $\delta$ basis result in $\Omega_2$:
$$\text{Mean Error}=\frac{\left|\sum_{i=1}^{P}\hat{\rho}_\Omega(\mathbf{r}_i)-\sum_{i=1}^{P}\rho_\Omega(\mathbf{r}_i)\right|}{\sum_{i=1}^{P}\rho_\Omega(\mathbf{r}_i)},$$
where $\hat{\rho}_\Omega$ is the approximation to the image source `exact' result $\rho_\Omega$. In configuration A the result of this mean error calculation is $8.3112\cdot 10^{-4}$
and so the approximation of the mean interior density in $\Omega_2$ is reasonably good (the mean error result for configuration B is $2.4259\cdot 10^{-3}$ and for configuration C is $1.4634\cdot 10^{-3}$). Furthermore, it is often only the mean response in a sub-domain that is of interest (depending on the application) and in statistical energy analysis, this is typically the only available quantity \cite{BM05, GT09}. The results in Fig.~\ref{fig:ConA} also suggest that our method will give reasonable results beyond this sub-domain average after applying local spatial averaging or low pass filtering to smooth the fluctuations. We expect the results to improve in the so-called ray chaos regime of longer trajectories where the dependence on the direction is milder and very large numbers of reflections must be included, which also means that the image source method soon becomes computationally prohibitive \cite{LS15}.

In order to explore the long-time ray behaviour further, we now investigate the effect of reducing the damping coefficient from $\mu=\pi/2$ to $\mu=\pi/200$. For this damping parameter, the image source method will still be used as a means of comparison. However, the image source results have now only converged to between 3 and 4 significant digits after 400 iterations and so we will use the term relative discrepancy as opposed to relative error to describe the differences between the methods. We note that convergence to 3 or 4 places will still be sufficient to identify an improvement in the $\mu=\pi/2$ relative error results described above.

\begin{figure}
\centering
\includegraphics[trim= 160  80 200 30, width=\textwidth]{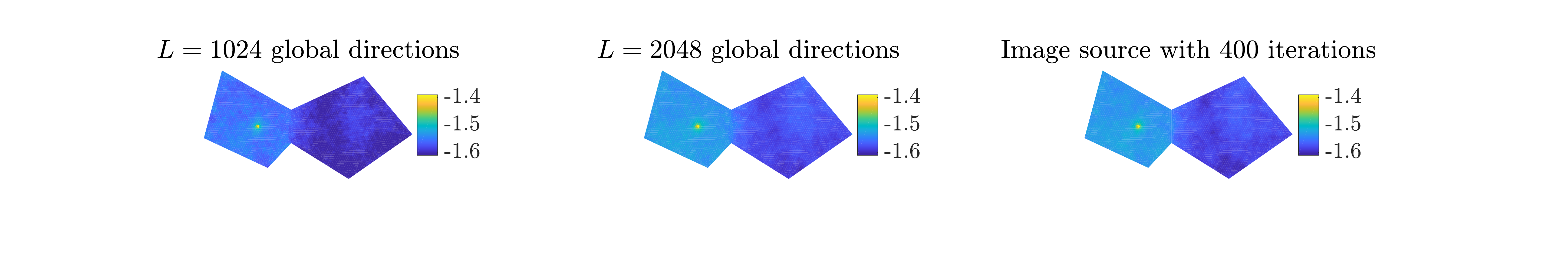}
\caption{Logarithm (base 10) of the interior energy densities in coupled cavity configuration A with damping coefficient $\mu=\pi/200$. The three sub-plots compare the results calculated using an average boundary element step size of $ 0.0125$ and $L=1024$ (left) or $L=2048$ (center) global ray directions with the result of applying 400 iterations of the image source method (right).}\label{fig:ConALD3}
\end{figure}

Figure \ref{fig:ConALD3} shows the results using the Dirac $\delta$ direction basis with both $L=1024$ and $L=2048$ global directions, and an average boundary element step size of $0.0125$, compared to the image source result described above. The match to the image source is clearly better with the increased number of directions and this is confirmed by the fact that the relative discrepancy (computed using the formula (\ref{MRE}) with the image source result in place of the exact solution) decreases from 0.04068 to 0.007755. We note that this is a considerable improvement on the results for larger damping and our earlier assertion about the accuracy of the method improving in the regime of many reflections appears to be correct. For completeness, the relative discrepancy for the more standard DEA approach with a Legendre polynomial basis approximation in the momentum variable up to degree 32 and an average boundary element step size of $0.0125$, as before, is 0.01985. The discrepancy is therefore in between the two results computed with the Dirac delta basis, but with far fewer degrees of freedom in the model. Note also that the computational overhead per degree of freedom for the Dirac delta basis is considerably less than for the Legendre polynomial basis since no numerical integration is necessary. \DC{The computation times for the results in Fig.~\ref{fig:ConALD3} again show an approximate $\mathcal{O}(L)$ scaling. A non-optimised serial MATLAB code gave computation times of 24 minutes with $L=1024$ and 49 minutes with $L=2048$, and as before, the time is dominated by the assembly cost for the linear system. All other tasks combined take less than 20s in total for both values of $L$ considered. } 

\subsection{Flexural vibrations of thin elastic shells}
\label{sec:flex-vibr-thin}
\subsubsection{Thin shell cylindrical ridge}
\begin{table}
\caption{\label{tab:ShellParams}\noindent Parameters for the thin elastic shell model.}
\centering
\begin{tabular}{|c|c|c|c|} \hline
Parameter description & Notation & Value & Units\\
\hline
 Material thickness &  $\zeta$ & 2.5e-3  & m \\
 Young's modulus &  $E$ & 7.0e10  & kg m$^{-1}$ s$^{-2}$\\
  Material density &  $\varrho$ & 2700  & kg m$^{-3}$\\
  Poisson's ratio  &  $\nu$ & 0.33  & - \\
    Angular frequency  &  $\omega$ &  18000$\pi$ & Rad\:s$^{-1}$\\
       \hline
  \end{tabular}
\end{table}
We first revisit the example of the curved ridge shown in Fig.~\ref{fig:Ridge0}, albeit with the underlying wave model replaced with the biharmonic equation for plate bending. Throughout this section the material parameters used correspond to a homogeneous thin aluminium shell, approximated locally as a flat plate in each triangular sub-domain, and are summarised in Table \ref{tab:ShellParams}. We consider a modified source density defined as 
\begin{equation}\label{eq:RidgeRho0}
\rho_0(s_j,p_j)=\mathbf{1}_{\{x_s=z_s=0\}}(s_j)\delta(p_j-p_0)
\end{equation}
for $j=1,2,\ldots,K$. The modified source term (\ref{eq:RidgeRho0}) corresponds to a line-source along the $y$-axis from $y=0$ to $y=1$ with tangential momentum $p_0=\eta\sin(\theta_0)$ and direction $\theta_0\in (-\pi/2,\pi/2)$. Note that the case $\theta_0=0$ corresponds to an incoming ray direction that is perpendicular to the boundary as considered in Section \ref{sec:simpeg}. We have dropped the subscript $j$ from $\eta$ since we consider a homogeneous aluminium shell and $$\eta=\frac{\sqrt{\omega}}{c}=\left(\frac{12\varrho(1-\nu^2)}{E\zeta^2}\right)^{1/4}$$ takes the value $\eta =0.5068$ (to four significant digits, see Table \ref{tab:ShellParams}) in every mesh cell $\Omega_j$, $j=1,2,\ldots, K$. Throughout this section we apply the simplest possible boundary element discretization on each mesh cell with exactly one boundary element per edge.

For flexural wave propagation on a thin shell cylindrical ridge, the bending wave only transmits across the ridge if the wavenumber in the direction of the maximum curvature is sufficiently large \cite{NS16}. In the example here, the initial direction of rays from the (upper-left) source edge needs to be below a critical value for transmission, otherwise the wave is reflected. For the ridge considered here this critical value can be computed from the known radius of curvature, taking the value $3\pi/2$ for $\tilde{x}\in(1/3,2/3)$  and zero elsewhere. When applying the technique on curved shells more generally, we instead numerically estimate the principle directions and curvatures for the shell using a nearest neighbour interpolation by a quadric surface \cite{TG06}.

Once the principle directions and associated curvatures are estimated, the threshold incoming ray direction for bending mode transmission may be calculated from the dispersion relationship for a flat plate and a cylindrical shell \cite{NS16}. In particular, the dispersion curve for a cylindrical shell is orthotropic and the second component of the wave vector $\mathbf{k}=(k_x, k_y)$ extends over a smaller range than in the dispersion curve for the flat region; waves with large $|k_y|$ in the flat region cannot pass into the curved region  (see Ref.~\cite{NS16} for details). The threshold incident ray directions can then be identified as corresponding to the maximum and minimum  admissible values of $k_y$ in the curved region. We note that by symmetry of the problem at hand only one of the extrema of $k_y$ must be calculated, since the other extremum will simply be the negative of this. The dispersion relation may be differentiated (analytically) to identify the extrema and then the threshold incoming direction \DC{$\theta^*$ at an interface between two mesh cells} may be calculated from the corresponding wave vector $\mathbf{k}=(k_x, k_y)$ via $\tan(\DC{\theta^*})=k_y/k_x$. \DC{Once the threshold directions have been calculated for each interface, they are implemented within $B_{I,J}$ (\ref{eq:TransferBij}) by setting $w_{i,j}(\tilde{p}_n(s_j))=1-\delta_{ij}$ for transmission from $\Omega_j$ to $\Omega_i$ and $w_{i,j}(\tilde{p}_n(s_j))=\delta_{ij}$ for reflection, where $\delta_{ij}$ is the Kronecker delta. If the incoming ray direction with respect to the principle direction associated to the maximum curvature is less than $\theta^*$ then we implement transmission. When the incoming ray direction exceeds $\theta^*$ then the direction of the ray reflected with respect to the principle direction is used to determine whether this ray corresponds to a reflection or transmission through the local interface. Note that the incoming ray direction with respect to the principle direction is calculated from $\tilde{p}_n(s_j)=\eta_j\sin(\phi_n(s_j))$ by associating $\phi_n$ with its corresponding global direction and using the scalar product of the global direction vector and the principle direction.}

\begin{figure}
\centering
\includegraphics[trim= 0 40  0  0, width=0.95\textwidth]{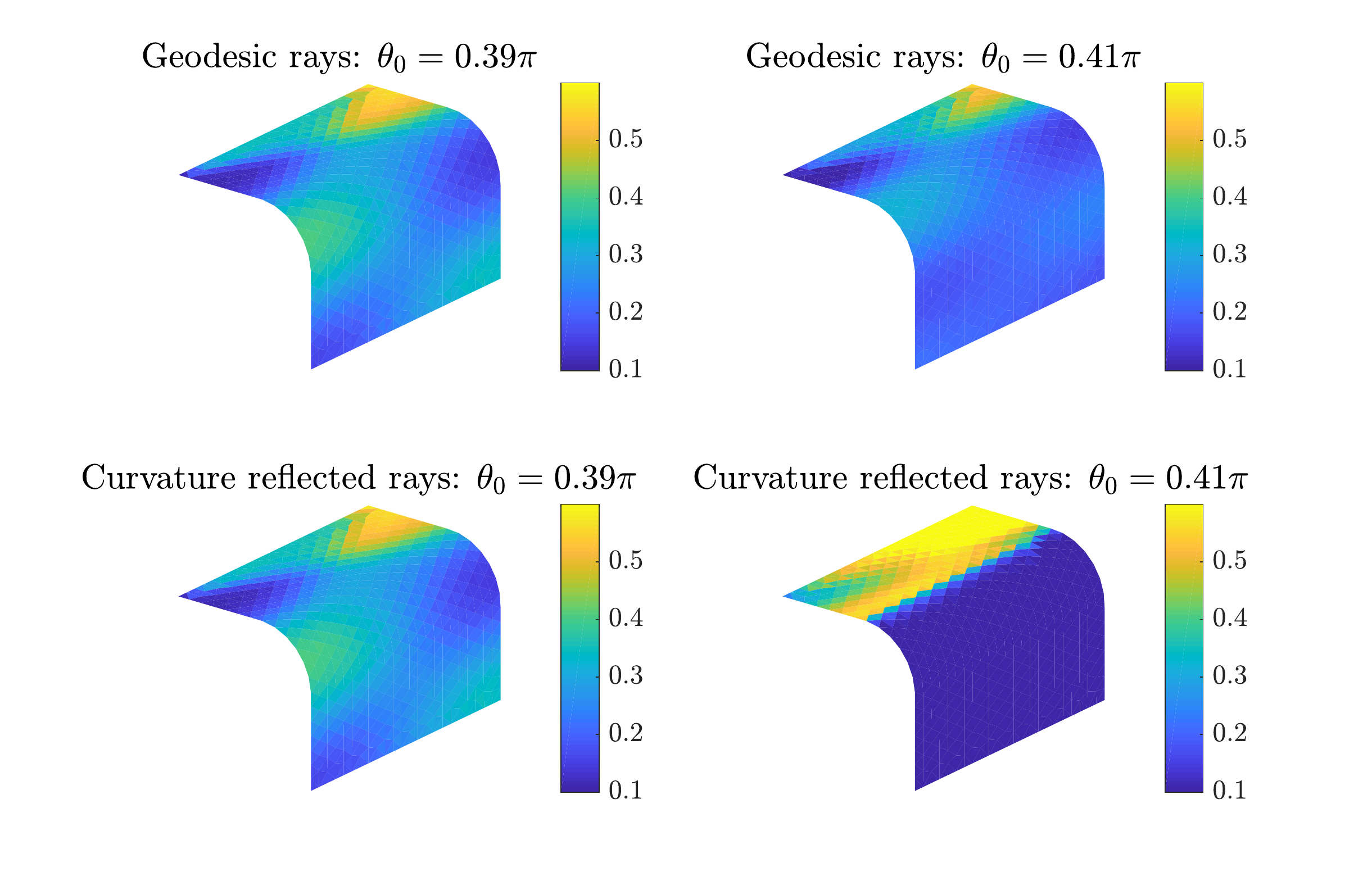}
\caption{\DC{Flexural wave energy densities on a cylindrical ridge approximated by $K=899$ triangular mesh cells and} with line sources emerging from the upper-left edge with direction $\theta_0=0.39\pi$ (left column) or $\theta_0=0.41\pi$ (right column). Top row: geodesic flow model with complete transmission across the ridge. Bottom row: model including reflections due to the curvature.}\label{fig:RidgeMulti}
\end{figure}

In the case of the aluminium ridge under consideration here with an incident wave of angular frequency $\omega=18000\pi$, the bending wavenumber is $k=120.527$ and the maximal $k_y$ value may be calculated as $114.698$, both to 6 significant digits. The threshold incident ray direction is therefore $\DC{\theta^*}=\tan^{-1}(k_y/\sqrt{k^2-k_y^2})\approx 0.4 \pi$. Figure \ref{fig:RidgeMulti} shows the results of incorporating these curvature induced reflections compared to a geodesic ray model with complete transmission across the ridge as considered previously. The choice of damping parameter $\mu$ here corresponds to $0.5\%$ hysteretic structural damping whereby $\mu=0.005k/2\approx 0.30132$. One observes that the results of both models are indistinguishable when $\theta_0=0.39\pi$ as would be expected since this choice of initial ray direction falls below the threshold value \DC{$\theta^*$}. However, for $\theta_0=0.41\pi$ the effect of the reflection model is evident since this value of $\theta_0$ exceeds the threshold value. Note that the global directions for the simulations in this section are all taken as $\Phi_l=2\pi (l-1)/L$, $l=1,2,...,L$ and for the calculations shown in Fig.~\ref{fig:RidgeMulti}, a choice of $L=200$ ensures that the directions corresponding to $\theta_0=0.39\pi$, $\theta_0=0.41\pi$ and all the possible boundary reflection directions are included within the basis. \DC{The computational cost for all four sub-plots is in the range 125s to 151s using the MATLAB parallel computing toolbox to run on 16 cores with a 2.6GHz clock speed. The times are again dominated by the linear system assembly (94s) and the longer times here correspond to the lower sub-plots, owing to a longer pre-processing time for calculating the local curvatures and principle directions.}

\subsubsection{Car shock tower}  We now apply the above described models to a thin aluminium shell with relevance to industrial applications, namely a shock tower model taken from a Range Rover car. This example was also considered in Ref.~\cite{CTLS13}, where the directional dependence of the ray density was approximated using a Legendre polynomial basis and the geometry induced reflections were modelled using plate junction theory, 
rather than curved shell theory as proposed here. The material properties are chosen to be the same as for the previous cylindrical ridge with constant thickness throughout. We note that this is a simplification of the true model, which  is instead comprised of several regions where the thickness of the shell differs between mesh cells \cite{CTLS13}. Including this thickness variation is an area for future work, which would also include devising a suitable direction basis for problems including refraction.

\begin{figure}
\centering
\includegraphics[trim= 90 40   80  0, width=\textwidth]{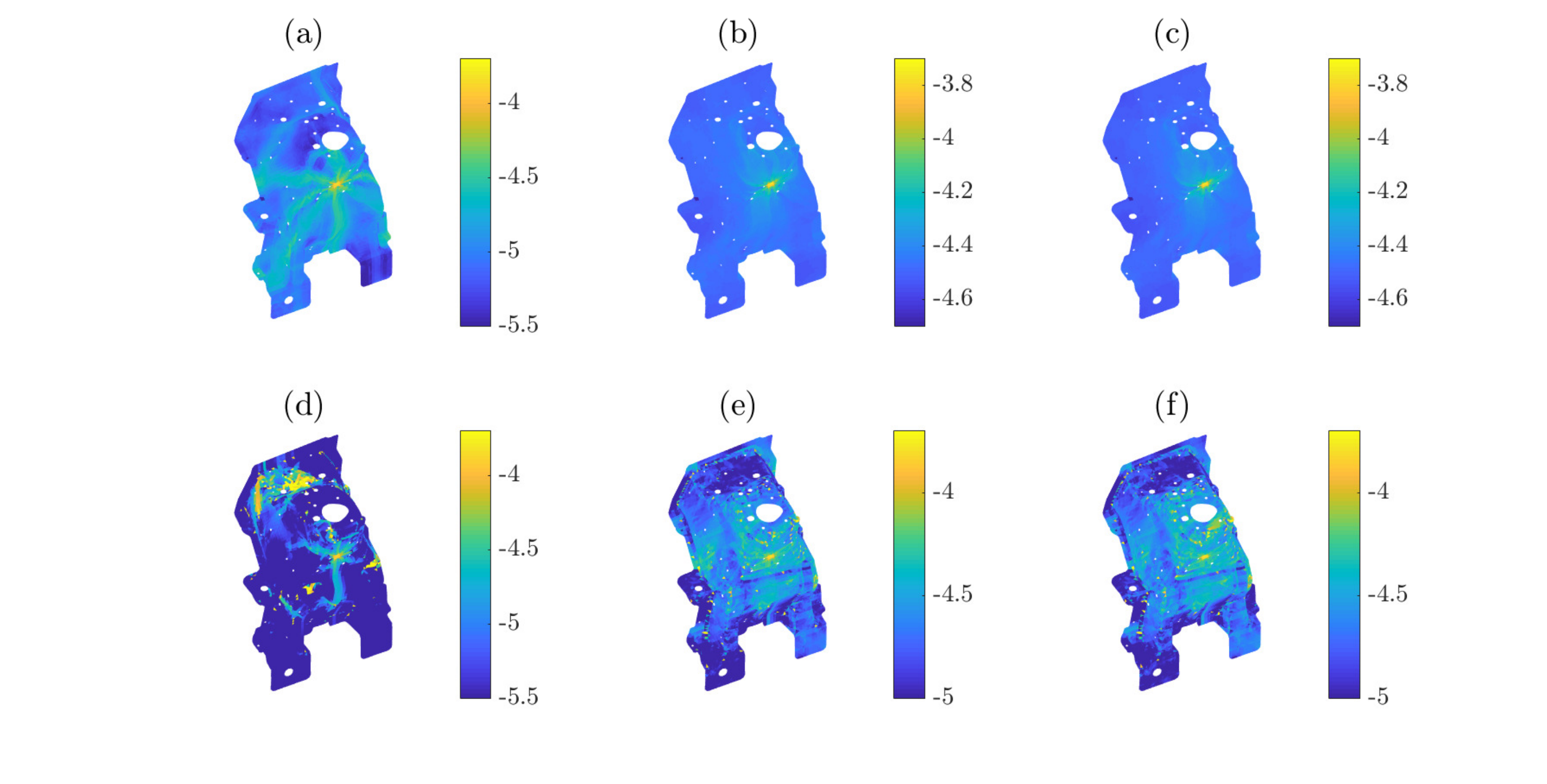}
\caption{Logarithm (base 10) of the interior energy densities for flexural wave motion excited by a point source  on a vehicle shock tower with damping coefficient $\mu=k/400$. The six sub-plots compare the results calculated using different numbers of ray directions in the direction basis as follows: (a) and (d) $L=8$, (b) $L=200$, (c) $L=400$, (e) \DC{$L=600$} and (f) \DC{$L=1000$}. Top row: geodesic flow model with complete transmission across the shell. Bottom row: model including reflections due to the curvature.}\label{fig:ST6RT}
\end{figure}

Figure \ref{fig:ST6RT} shows the results for the approximated ray density distribution (\ref{eq:RhoOmegaDS})  calculated at the centroids of each of the $K=11623$ triangular mesh cells used to approximate the geometry. The excitation of the system is via a unit point source that is visible on the front face of the shock tower and takes the form \cite{TH19}:
\begin{equation}\label{eq:Rho0ptv}
\rho_0(s_j,p_j;\mathbf{r}_0)=\frac{ k^2 \cos(\theta_0)e^{-\mu D(\mathbf{r}_0,s_j)}w_{j,0}(s_j)\delta(p_j-p_0)}{16\pi\zeta\varrho\omega^{3/2} D(\mathbf{r}_0,s_j)},
\end{equation}
where $\mathbf{r}_0\in\mathbb{R}^3$ contains the Cartesian coordinates of the source location. In this case the source point corresponds to a vertex of the mesh and so the weight function $w_{j,0}$ is essentially a visibility function for the mesh cells $\Omega_j$ whose boundaries can receive rays directly from $\mathbf{r}_0$. In particular, assuming that $\mathbf{r}_0$ is not on a vertex connecting to a free edge and that the region surrounding to the source point is flat and homogeneous then $w_{j,0}(s_j)=0$ unless either $\mathbf{r}_0\in\Gamma_j$ and $s_j$ is on a free edge or $\mathbf{r}_0\not\in\Gamma_j$, but $s_j$ is on  an edge shared with a mesh cell whose boundary contains $\mathbf{r}_0$. In these latter two cases, $w_{j,0}(s_j)=1$. The choice of damping parameter $\mu$ again corresponds to $0.5\%$ hysteretic structural damping with $\mu\approx 0.30132$.

Note that unlike the cylindrical ridge considered in the previous section, here there is no simple and consistent choice of local $x$-axis that can be applied in each mesh cell to correspond to the global direction $\Phi_1 =0$. Recall that for the ridge example, the curvilinear coordinate $\tilde{x}$ takes the role of the $x$-coordinate and so the global directions in each mesh cell can be chosen relative to $\tilde{x}$ (after projection onto the mesh cell under consideration). Here we employ a more generally applicable procedure that preserves the propagation directions on flat regions of the structure, that is, for coplanar mesh cells. To this end we perform a rotation of the global Cartesian coordinate frame (in $\mathbb{R}^3$) so that the rotated $z$-axis always corresponds to a consistently defined normal vector across the whole surface $\Omega$ (the angle of rotation for a given mesh cell is therefore the angle between the normal vector to that mesh cell and the positive $z$-axis). Uniqueness of the rotation follows by defining the axis of rotation to be the line formed from the intersection of the tangent plane to the mesh cell with the plane $z=0$. The global direction $\Phi_1 =0$ then corresponds to the rotated $x$-axis and the term global direction becomes a slight misnomer due to the local rotations for non-coplanar mesh cells.

The top row of Fig.~\ref{fig:ST6RT} shows the results of a geodesic flow model with complete transmission across all mesh edges that do not represent a physical edge or hole in the structure. One can gain some intuition from plot (a), which shows the model using only $L=8$ global directions in the direction basis. These 8 directions can be seen emanating from the source point on the front face of the shock tower and some of the early reflections from holes and the outer boundary can also be observed. Plots (b) and (c) show the same model when the number of global directions is increased to $L=200$ and $L=400$, respectively. The latter two plots appear identical by eye demonstrating convergence of the method and show that the vibrational energy spreads evenly across the shock tower from the source point. The energy is noticeably higher in the vicinity of the source point due to the damping. The bottom row of Fig.~\ref{fig:ST6RT} shows the results of additionally including reflections due to the curvature in the model. Plot (d) may be directly compared with plot (a) since only $L=8$ global directions are modelled and the colour scales are consistent. One observes that, as for the cylindrical ridge, including curvature induced reflections leads to the localisation of the vibrational energy on relatively flat regions of the structure. Plots (e) and (f) show the same model when the number of global directions is increased to \DC{$L=600$} and \DC{$L=1000$}, respectively. The latter two plots appear very similar in their structure, \DC{providing evidence of convergence}. The inclusion of curvature induced reflections in the model therefore significantly increases the number of global directions that are necessary for convergence to be observed in the simulations. The results in plots (e) and (f) demonstrate that vibrational energy localises on relatively flat regions of the shell and the curved ridges appear to suppress vibrations. This is consistent with observations from the finite element simulations of the shock tower in Ref.~\cite{CTLS13}.

\DC{The total computational cost for the results here ranges from 9 hours ($L=8$) to 2.5 days ($L=1000$)  for non-optimised parallel MATLAB code running on 16 cores as before. These times remain dominated by the linear system assembly, which ranges from 6.7 hours (for both $L=8$ and $L=200$) to 34 hours (for $L=1000$). The linear system solver now also takes a significant proportion of the computation time and so we switched from the direct solver used previously (MATLAB backslash) to an iterative solver (\texttt{bicgstab} as employed for large DEA problems in Refs.~\cite{CTLS13, TH19}). The linear solver times range  between approximately 1 hour ($L=8$) and 1 day ($L=1000$).}

\section{Conclusions}

A Petrov-Galerkin discretization of a phase-space boundary integral equation for transporting wave energy densities has been presented that has the property of preserving the direction of propagation across multi-component domains in $\mathbb{R}^2$. The property relies on the existence of a consistent global direction set and extends to curved surfaces when the geodesics can be used as a reference for the global direction set. More generally, if the geometry is complex and the geodesics are not known \emph{a-priori} then one can define a transformation of the global direction set into each mesh cell that can, at least, maintain the direction preservation property within the flat regions of the shell. Numerical experiments have demonstrated that the direction preservation property is important for accurately solving problems with relatively simple ray dynamics, as well as for correctly capturing direction dependent reflections due to curvature on thin shells. The global direction set for this work was chosen to be equally spaced with rotational symmetry in each quartile, but we note that there is freedom to choose these directions according to any problem specific knowledge that you may have; for example, a geometry with some form of symmetry or a structure with waveguides in particular directions. The method was also shown to produce accurate results in the so-called ray chaos regime where the directional dependence of the solution becomes smoothed. Finally, an example from industry was considered, demonstrating the effectiveness of the method for simulating vibrational energy distributions in geometrically complex shell structures. 

One of the main advantages of the method proposed here is that the discretization of the governing linear integral operator can be done without numerical integration; the matrix entries for the associated linear system are computed using a single analytically calculated integral. These relatively simple calculations can be performed in parallel many times, meaning that large and industrially relevant examples are within reach. Comparing the methodology proposed here with existing DEA approaches, we see that in this work we are able to accurately solve the highly directive problems that are typically problematic for DEA with relatively small numerical models. We are also able to accurately simulate the ray chaotic situations that existing DEA techniques are able to model most efficiently, but with larger numerical models including many ray directions. \DC{However, the computational cost of the proposed method scales only linearly with the number of ray directions in the direction basis owing to the sparsity of the projection onto this basis. This compares favourably with the quadratic (or worse) scaling with the maximum degree polynomial in the direction basis for previous DEA approaches. The development of hybrid schemes including both orthogonal polynomials and delta distributions in the approximation space to combine the advantages of both approaches could therefore be a promising area for future work.} 
\appendix
\section{From waves to rays}\label{sec:App}
We consider frequency domain wave equations of the form:
\begin{equation}\label{eq:Wave}
\Delta^{\alpha} u - (\mathrm{i}k)^{2\alpha} u =0,
\end{equation}
in the large angular frequency regime $\omega\gg 1$, where $k=\omega/c$ is the wavenumber with $c>0$, the wave speed. When $\alpha=1$ and $c = \mathrm{const}$, the wave equation (\ref{eq:Wave}) corresponds to the Helmholtz equation. If instead $\alpha=2$ and $c(k) \propto k$, then (\ref{eq:Wave}) is the biharmonic wave equation for modelling flexural plate motion. The dispersion property $\omega\propto k^\alpha$ leads to seeking solutions in the form
$$u(\mathbf{r})=e^{\mathrm{i}\omega^{(1/\alpha)}S(\mathbf{r})}\sum_{\kappa=0}^\infty A_\kappa(\mathbf{r}) (\mathrm{i}\omega)^{-\kappa}$$
for large $\omega$. Substitution into the wave equation (\ref{eq:Wave}) results in the Eikonal equation at leading order $\mathcal{O}(\omega^2)$ for the phase function $S$:
$$|\nabla S|=\frac{\omega^{(\alpha-1)/2}}{c}=:\eta.$$
One then obtains a Hamiltonian system of ODEs (ray equations) for the trajectories $Y=(\mathbf{r},\mathbf{p})=(\mathbf{r},\nabla S)$ associated to the Hamiltonian
$$H(\mathbf{r},\mathbf{p})=|\mathbf{p}|/\eta(\mathbf{r})=1$$
via the method of characteristics.

For frequency domain wave problems governed by (\ref{eq:Wave}), it is the long-time evolution of the Hamiltonian system that is of interest, which is often studied in terms of the associated (particle) density distribution $f$ in phase-space. The density distribution obeys a phase-space conservation law known as the Liouville equation:
\begin{equation}\label{eq:Liouville}
\frac{\partial f}{\partial t}(Y,t)+\frac{\partial
Y}{\partial t}\cdot\nabla_{Y}(f(Y,t))=0,
\end{equation}
and the long-time behaviour may be conveniently studied in terms of the associated stationary problem with $\partial f/\partial t\equiv 0$. The method of characteristics yields an expression for the solution of (\ref{eq:Liouville}) in terms of the Frobenius-Perron operator (\ref{eq:FPO}), and consequently the stationary problem too by considering the density accumulated from a continuous source in the long-time limit (\ref{eq:FPOstat}).
\bibliographystyle{siam}
\bibliography{deltaDEA}

\begin{thebibliography}{10}

\bibitem{JB17}
{\sc J.~Bajars, D.J. Chappell, T.~Hartmann, and G.~Tanner}, {\em Improved
  approximation of phase-space densities on triangulated domains using discrete
  flow mapping with p-refinement}, J.~Sci~Comput, 72 (2017), pp.~1290--1312.

\bibitem{Lebot2010}
{\sc A.~Le Bot and V.~Cotoni}, {\em Validity diagrams of statistical energy
  analysis}, J.~Sound Vib., 329 (2010), pp.~221 -- 235.

\bibitem{BMM12}
{\sc M.~Budi\u{s}i\'{c}, R.~Mohr, and I.~Mezi\'{c}}, {\em Applied
  {K}oopmanism}, Chaos, 22 (2012), p.~047510.

\bibitem{CTLS13}
{\sc D.J. Chappell, G.~Tanner, D.~L{\"o}chel, and N.~S{\o}ndergaard}, {\em
  Discrete flow mapping: transport of phase space densities on triangulated
  surfaces}, P.~Roy.~Soc.~Lond.~A~Mat, 469 (2013).

\bibitem{JC09}
{\sc J.~Cuenca, F.~Gautier, and Simon L.}, {\em The image source method for
  calculating the vibrations of simply supported convex polygonal plates},
  J.~Sound Vib., 322 (2009), pp.~1048--1069.

\bibitem{PD11}
{\sc P.~Dutta and R.~Bhattacharya}, {\em Hypersonic state optimisation using
  the {F}robenius-{P}erron operator}, J.~Guid.~Control Dynam., 34 (2011),
  pp.~325--333.

\bibitem{BE03}
{\sc B.~Engquist and O.~Runborg}, {\em Computational high frequency wave
  propagation}, Acta Numer., 12 (2003), pp.~181--266.

\bibitem{GF09}
{\sc G.~Friesecke, O.~Junge, and P.~Koltai}, {\em Mean field approximation in
  conformation dynamics}, Multiscale Model.~Simul., 8 (2009), pp.~254--268.

\bibitem{FJK11}
{\sc G.~Froyland, O.~Junge, and P.~Koltai}, {\em Estimating long term behavior
  of flows without trajectory integration: the infinitesimal generator
  approach}, SIAM J. Numer. Anal., 51 (2013), pp.~223--247.

\bibitem{TG06}
{\sc T.~Gatzke and C.M. Grimm}, {\em Estimating curvature on triangular
  meshes}, Int.~J. Shape Modeling, 12 (2006), pp.~1 -- 28.

\bibitem{Gradoni2014}
{\sc G.~Gradoni, J.~Yeh, B.~Xiao, T.M. Antonsen, S.M. Anlage, and E.~Ott}, {\em
  Predicting the statistics of wave transport through chaotic cavities by the
  random coupling model: A review and recent progress}, Wave Motion, 51 (2014),
  pp.~606 -- 621.

\bibitem{TH19}
{\sc T.~Hartmann, S.~Morita, G.~Tanner, and D.J. Chappell}, {\em High-frequency
  structure- and air-borne sound transmission for a tractor model using
  dynamical energy analysis}, Wave Motion, 87 (2019), pp.~132 -- 150.
\newblock Innovations in Wave Modelling II.

\bibitem{QH19}
{\sc Q.~Huang, J.~Ding, and N.H. Rhee}, {\em Convergence analysis of projection
  methods for {F}robenius-{P}erron operators based on matrix norm techniques},
  J.~Math.~Anal.~Appl., 479 (2019), pp.~337 -- 349.

\bibitem{JK09}
{\sc O.~Junge and P.~Koltai}, {\em Discretization of the {F}robenius-{P}erron
  operator using a sparse haar tensor basis: The sparse {U}lam method}, SIAM J.
  Numer. Anal., 47 (2009), pp.~3464--3485.

\bibitem{KKS16}
{\sc S.~Klus, P.~Koltai, and C.~Sch{\"u}tte}, {\em On the numerical
  approximation of the {P}erron-{F}robenius and {K}oopman operator},
  J.~Comp.~Dyn., 3 (2016), pp.~57--79.

\bibitem{KKS18}
{\sc S.~Klus, F.~N{\"u}ske, P.~Koltai, H.~Wu, I.~Kevrekidis, C.~Sch{\"u}tte,
  and F.~No{\'e}}, {\em Data-driven model reduction and transfer operator
  approximation}, J.~Nonlinear Sci., 28 (2018), pp.~985--1010.

\bibitem{JK16}
{\sc J.~Kullig and J.~Wiersig}, {\em Frobenius-{P}erron eigenstates in deformed
  microdisk cavities: non-{H}ermitian physics and asymmetric backscattering in
  ray dynamics}, New J.~Phys., 18 (2016), p.~015005.

\bibitem{HK09}
{\sc H.~Kuttruff}, {\em Room acoustics}, Spon Press, UK, 5th~ed., 2009.

\bibitem{Lafont2013}
{\sc T.~Lafont, N.~Totaro, and A.~Le~Bot}, {\em Review of statistical energy
  analysis hypotheses in vibroacoustics}, Proc. R. Soc. A, 470 (2013).

\bibitem{HL89}
{\sc H.~Ling, R.~Chou, and S.~Lee}, {\em Shooting and bouncing rays:
  calculating the rcs of an arbitrarily shaped cavity}, IEEE Transactions on
  Antennas and Propagation, 37 (1989), pp.~194--205.

\bibitem{WL16}
{\sc W.~Lu, J.~Qian, and R.~Burridge}, {\em Babich's expansion and the fast
  {H}uygens sweeping method for the {H}elmholtz equation at high frequencies},
  J.~Comp.~Phys., 313 (2016), pp.~478 -- 510.

\bibitem{Lyon1969}
{\sc R.H. Lyon}, {\em Statistical analysis of power injection and response in
  structures and rooms}, J.~Acoust.~Soc.~Am., 45 (1969), pp.~545--565.

\bibitem{RL95}
{\sc R.H. Lyon and R.G. DeJong}, {\em Theory and Application of Statistical
  Energy Analysis}, Butterworth-Heinemann, Boston, 1995.

\bibitem{BM05}
{\sc B.R. Mace}, {\em Statistical energy analysis: coupling loss factors,
  indirect coupling and system modes}, J.~Sound Vib., 279 (2005), pp.~141 --
  170.

\bibitem{RN18}
{\sc R.A. Norton, C.~Fox, and M.E. Morrison}, {\em Numerical approximation of
  the {F}robenius-{P}erron operator using the finite volume method}, SIAM
  J.~Numer.~Anal., 56 (2018), pp.~570--589.

\bibitem{distmesh}
{\sc P.O. Persson and G.~Strang}, {\em A simple mesh generator in {MATLAB}},
  SIAM Rev., 46 (2004), pp.~329--345.

\bibitem{JQ16}
{\sc J.~Qian, W.~Lu, L.~Yuan, S.~Luo, and R.~Burridge}, {\em Eulerian
  geometrical optics and fast huygens sweeping methods for three-dimensional
  time-harmonic high-frequency maxwell's equations in inhomogeneous media},
  Multiscale Model.~Sim., 14 (2016), pp.~595--636.

\bibitem{FR19}
{\sc F.~Rullan and M.M. Betcke}, {\em Hamilton-{G}reen solver for the forward
  and adjoint problems in photoacoustic tomography}, Arxiv preprint,
  arXiv:1810.13196v1 (2018).

\bibitem{Froyland10}
{\sc N.~Santitissadeekorn, G.~Froyland, and A.~Monahan}, {\em Optimally
  coherent sets in geophysical flows: A transfer-operator approach to
  delimiting the stratospheric polar vortex}, Phys.~Rev.~E, 82 (2010),
  p.~056311.

\bibitem{LS15}
{\sc L.~Savioja and U.P. Svensson}, {\em Overview of geometrical room acoustic
  modeling techniques}, J.~Acoust.~Soc.~Am., 138 (2015), pp.~708--730.

\bibitem{JS2020}
{\sc J.~Slipantschuk, M.~Richter, D.J. Chappell, and G.~Tanner}, {\em Transfer
  operator approach to ray-tracing in circular domains}, Nonlinearity, 33
  (2020), pp.~5773--5790.

\bibitem{IS01}
{\sc I.P. Smirnov, A.L. Virovlyansky, and G.M. Zaslavsky}, {\em Theory and
  applications of ray chaos to underwater acoustics}, Phys. Rev. E, 64 (2001),
  p.~036221.

\bibitem{NS16}
{\sc N.~Sondergaard and D.J. Chappell}, {\em Ray and wave scattering in
  smoothly curved thin shell cylindrical ridges}, J.~Sound Vib., 377 (2016),
  pp.~155 -- 168.

\bibitem{GT09}
{\sc G.~Tanner}, {\em Dynamical energy analysis -- determining wave energy
  distributions in vibro-acoustical structures in the high-frequency regime},
  J.~Sound Vib., 320 (2009), pp.~1023 -- 1038.

\bibitem{Ulam1964}
{\sc S.M. Ulam}, {\em Problems in Modern Mathematics}, John Wiley and Sons, New
  York, 1964.

\bibitem{Cer01}
{\sc V.~\v{C}erven\'y}, {\em Seismic Ray Theory}, Cambridge University Press,
  Cambridge, 2001.

\bibitem{CW19}
{\sc C.~Wormell}, {\em Spectral {G}alerkin methods for transfer operators in
  uniformly expanding dynamics}, Numer.~Math., 142 (2019), pp.~421--463.

\end{thebibliography}
\end{document}